\documentclass[11pt]{article}
\usepackage[mathscr]{eucal}
\usepackage{mathrsfs}
\usepackage{amssymb}
\usepackage{calligra}
\usepackage{fontenc}
\usepackage{amsbsy} 
\usepackage{graphicx}
\graphicspath{%
    {converted_graphics/}
    {/}
}
\begin{document}

\newcommand{\ts}{{\tilde{\sf s}}}
\newcommand{\sfv}{{\sf v}}
\newcommand{\sfw}{{\sf w}}
\newcommand{\simge}{\ba{cc}\vspace*{-2.4mm}>\\ \sim\ea }
\newcommand{\simle}{\ba{cc}\vspace*{-2.4mm}<\\ \sim\ea }
\newcommand{\Cdot}{\!\cdot\!}
\newcommand{\sq}{{$\sqcap\!\!\!\!\sqcup$}}
\newcommand{\Eu}{{\rm I\,\!\! E}}
\newcommand{\Io}{\Int{\Omega}{}}
\newcommand{\Id}{\Int{\cald}{}}
\newcommand{\Div}{\mbox{\rm div}\,}
\newcommand{\tr}{\mbox{\rm tr}\,}
\newcommand{\grad}{\mbox{\rm grad}\,}
\newcommand{\supp}{\mbox{\rm supp}\,}
\newcommand{\curl}{\mbox{\rm curl}\,}
\newcommand{\Ido}{\Int{\partial\Omega}{}}
\newcommand{\IdS}{\Int{\Sigma}{}}
\newcommand{\Oint}[2]{{\displaystyle \oint_{#1}^{#2}}}
\newcommand{\Int}[2]{{\displaystyle \int_{ #1}^{ #2}}}
\newcommand{\Lim}[1]{{\displaystyle \lim_{ #1}}}
\newcommand{\Limsup}[1]{{\displaystyle \limsup_{\footnotesize #1}}}
\newcommand{\Liminf}[1]{{\displaystyle \liminf_{\footnotesize #1}}}
\newcommand{\Sup}[1]{{\displaystyle \sup_{#1}}}
\newcommand{\Inf}[1]{{\displaystyle \inf_{#1}}}
\newcommand{\Max}[1]{{\displaystyle \max_{#1}}}
\newcommand{\Min}[1]{{\displaystyle \min_{#1}}}
\newcommand{\Sum}[2]{{\displaystyle \sum_{#1}^{#2}}}
\newcommand{\Prod}[2]{{\displaystyle \prod_{#1}^{#2}}}
\newcommand{\BCup}[2]{{\displaystyle \bigcup_{#1}^{#2}}}
\newcommand{\BCap}[2]{{\displaystyle \bigcap_{#1}^{#2}}}
\newcommand{\Frac}[2]{\displaystyle{\frac{\displaystyle{#1}}{\displaystyle{#2}}}}
\newcommand{\norm}[1]{\left\|{#1}\right\|}
\newcommand{\Norm}[1]{\langle\langle{#1}\rangle\rangle_q}
\newcommand{\No}[1]{\langle\!\langle{#1}\rangle\!\rangle}
\newcommand{\NO}[1]{{\langle{#1}\rangle}_{\lambda,q}}
\newcommand{\beea}{\begin{eqnarray}}
\newcommand{\eeea}{\end{eqnarray}}
\newcommand{\ms}{\medskip\smallskip}
\newcommand{\bs}{\bigskip}
\newcommand{\ps}{\par\smallskip}
\newcommand{\bfe}{{\mbox{\boldmath $e$}} }
\newcommand{\pni}{{\par\noindent}}
\newcommand{\bfq}{{\mbox{\boldmath $q$}} }
\newcommand{\bfz}{{\mbox{\boldmath $z$}} }
\newcommand{\0}{{\mbox{\boldmath $0$}} }
\newcommand{\LE}{\!\!\!&\le&\!\!\!}
\newcommand{\BL}[1]{{\par\smallskip{\bf Lemma #1.}}}
\newcommand{\BT}[1]{{\par\smallskip{\bf Theorem #1.}}}
\newcommand{\Ln}{[\!|}
\newcommand{\Rn}{|\!]}
\newcommand{\n}[1]{{\Ln{#1}\Rn}} 
\newcommand{\nq}[1]{{\Ln{#1}\Rn}_{q}} 
\newcommand{\nqr}[1]{{\Ln{#1}\Rn}_{q,r}} 
\newcommand{\Nq}[1]{{\langle{#1}\rangle}_{q}} 
\newcommand{\Nql}[1]{{\langle{#1}\rangle}_{\lambda,q}} 
\newcommand{\Nqr}[1]{{\langle{#1}\rangle}_{q,r}}
\newcommand{\N}[1]{{|\!\!|\!\!|\,{#1}\,|\!|\!\!|_2}}
\newcommand{\EA}[2]{$$#1$$%
\vspace{-6.mm}
\begin{equation}
\end{equation}
\vspace{-6.mm}
$$
#2
\setlength{\belowdisplayskip}{3mm}
\setlength{\belowdisplayshortskip}{3mm}
$$
}
\newcommand{\A}[2]{$$#1$$%
\vspace{-4.mm}
$$
#2
\setlength{\belowdisplayskip}{3mm}
\setlength{\belowdisplayshortskip}{3mm}
$$
}
\newcommand{\BF}{\begin{footnotesize}}
\newcommand{\EF}{\end{footnotesize}}
\setlength{\jot}{.15in}
\newcommand{\pde}[2]{{\displaystyle \frac{\mbox{$\partial #1$}}{\mbox{$\partial #2$}}}}
\newcommand{\ode}[2]{{\displaystyle \frac{\mbox{$d #1$}}{\mbox{$d #2$}}}}
\newcommand{\f}[2]{\frac{\mbox{$#1$}}{\mbox{$ #2$}}}
\newcommand{\bi}{\begin{itemize}}
\newcommand{\ei}{\end{itemize}}
\newcommand{\ed}{\end{document}}
\newcommand{\be}{\begin{equation}}
\newcommand{\ba}{\begin{array}}
\newcommand{\ea}{\end{array}}
\newcommand{\ee}{\end{equation}}
\newcommand{\eeq}[1]{\label{eq:#1}\end{equation}}
\newcommand{\real}{{\mathbb R}}
\newcommand{\compl}{{\mathbb C}}
\def\Id{\mbox{\boldmath $1$}}
\def\zero{\mbox{\boldmath $0$}}
\newcommand{\PP}{{\rm I\!\!\,P}}
\newcommand{\nat}{{\mathbb N}}
\newcommand{\bfpsi}{\mbox{\boldmath $\psi$}}
\newcommand{\bfchi}{\mbox{\boldmath $\chi$}}
\newcommand{\bfomega}{\mbox{\boldmath $\omega$}}
\newcommand{\bfome}{\mbox{\boldmath $\varpi$}}
\newcommand{\bfvaromega}{\mbox{\boldmath $\varpi$}}
\newcommand{\bfOmega}{\mbox{\boldmath $\Omega$}}
\newcommand{\bfTheta}{\mbox{\boldmath $\Theta$}}
\newcommand{\bfxi}{\mbox{\boldmath $\xi$}}
\newcommand{\bfmu}{\mbox{\boldmath $\mu$}}
\newcommand{\bfx}{\mbox{\boldmath $x$}}
\newcommand{\bfy}{\mbox{\boldmath $y$}}
\newcommand{\bfPsi}{\mbox{\boldmath $\Psi$}}
\newcommand{\bfphi}{\mbox{\boldmath $\varphi$}}
\newcommand{\bfhi}{\mbox{\boldmath $\phi$}}
\newcommand{\bfPhi}{\mbox{\boldmath $\Phi$}}
\newcommand{\bfv}{{\mbox{\boldmath $v$}} }
\newcommand{\bfu}{{\mbox{\boldmath $u$}} }
\newcommand{\bfsf}{{\mbox{\footnotesize\boldmath $s$}} }
\newcommand{\bfuf}{{\mbox{\footnotesize\boldmath $u$}} }
\newcommand{\bfw}{{\mbox{\boldmath $w$}} }
\newcommand{\bff}{{\mbox{\boldmath $f$}} }
\newcommand{\bfa}{{\mbox{\boldmath $a$}} }
\newcommand{\bfi}{{\mbox{\boldmath $i$}} }
\newcommand{\bfj}{{\mbox{\boldmath $j$}} }
\newcommand{\bfc}{{\mbox{\boldmath $c$}} }
\newcommand{\bfo}{{\mbox{\boldmath $o$}} }
\newcommand{\bfp}{{\mbox{\boldmath $p$}} }
\newcommand{\bfkp}{{\mbox{\footnotesize{\boldmath $k$}}} }
\newcommand{\bfka}{{\mbox{\footnotesize{\boldmath $k^*$}}} }
\newcommand{\bft}{{\mbox{\boldmath $t$}} }
\newcommand{\bfd}{{\mbox{\boldmath $d$}} }
\newcommand{\bfl}{{\mbox{\boldmath $l$}} }
\newcommand{\bfr}{{\mbox{\boldmath $r$}} }
\newcommand{\bfk}{{\mbox{\boldmath $k$}} }
\newcommand{\bfA}{{\mbox{\boldmath $A$}} }
\newcommand{\bfS}{{\mbox{\boldmath $S$}} }
\newcommand{\bfO}{{\mbox{\boldmath $O$}} }
\newcommand{\bfM}{{\mbox{\boldmath $M$}} }
\newcommand{\bfP}{{\mbox{\boldmath $P$}} }
\newcommand{\bfB}{{\mbox{\boldmath $B$}} }
\newcommand{\bfR}{{\mbox{\boldmath $R$}} }
\newcommand{\bfC}{{\mbox{\boldmath $C$}} }
\newcommand{\bfD}{{\mbox{\boldmath $D$}} }
\newcommand{\bfQ}{{\mbox{\boldmath $Q$}} }
\newcommand{\bfZ}{{\mbox{\boldmath $Z$}} }
\newcommand{\bfG}{{\mbox{\boldmath $G$}} }
\newcommand{\bfE}{{\mbox{\boldmath $E$}} }
\newcommand{\bfX}{{\mbox{\boldmath $X$}} }
\newcommand{\bfY}{{\mbox{\boldmath $Y$}} }
\newcommand{\bfH}{{\mbox{\boldmath $H$}} }
\newcommand{\bfI}{{\mbox{\boldmath $I$}} }
\newcommand{\bfJ}{{\mbox{\boldmath $J$}} }
\newcommand{\bfN}{{\mbox{\boldmath $N$}} }
\newcommand{\bfh}{{\mbox{\boldmath $h$}} }
\newcommand{\bfm}{{\mbox{\boldmath $m$}} }
\newcommand{\bfone}{{\mbox{\boldmath $1$}} }
\newcommand{\hs}{{\rm I}\!\!\,{\rm R}^3_+}
\newcommand{\cala}{{\cal A}}
\newcommand{\calb}{{\cal B}}
\newcommand{\calc}{{\cal C}}
\newcommand{\cald}{{\cal D}}
\newcommand{\cale}{{\cal E}}
\newcommand{\calf}{{\cal F}}
\newcommand{\calg}{{\cal G}}
\newcommand{\calh}{{\cal H}}
\newcommand{\cali}{{\cal I}}
\newcommand{\calj}{{\cal J}}
\newcommand{\calk}{{\cal K}}
\newcommand{\call}{{\cal L}}
\newcommand{\calm}{{\cal M}}
\newcommand{\caln}{{\cal N}}
\newcommand{\calo}{{\cal O}}
\newcommand{\calp}{{\cal P}}
\newcommand{\calq}{{\cal Q}}
\newcommand{\calr}{{\cal R}}
\newcommand{\cals}{{\cal S}}
\newcommand{\calt}{{\cal T}}
\newcommand{\calu}{{\cal U}}
\newcommand{\calv}{{\cal V}}
\newcommand{\calx}{{\cal X}}
\newcommand{\caly}{{\cal Y}}
\newcommand{\calw}{{\cal W}}
\newcommand{\calz}{{\cal Z}}
\newcommand{\bfsigma}{\mbox{\boldmath $\sigma$}}
\newcommand{\bfSigma}{\mbox{\boldmath $\Sigma$}}
\newcommand{\bftau}{\mbox{\boldmath $\tau$}}
\newcommand{\bfeta}{\mbox{\boldmath $\eta$}}
\newcommand{\bfT}{{\mbox{\boldmath $T$}} }
\newcommand{\bfV}{{\mbox{\boldmath $V$}} }
\newcommand{\bfU}{{\mbox{\boldmath $U$}} }
\newcommand{\bfW}{{\mbox{\boldmath $W$}} }
\newcommand{\bfF}{{\mbox{\boldmath $F$}} }
\newcommand{\bfK}{{\mbox{\boldmath $K$}} }
\newcommand{\bfL}{{\mbox{\boldmath $L$}} }
\newcommand{\bfb}{{\mbox{\boldmath $b$}} }
\newcommand{\bfg}{{\mbox{\boldmath $g$}} }
\newcommand{\bfn}{{\mbox{\boldmath $n$}} }
\newcommand{\bfs}{{\mbox{\boldmath $s$}} }
\newcommand{\cf}{{\it cf.} }
\newcommand{\io}{\int_\Omega}
\newcommand{\1}{\item[({\it i})]}
\newcommand{\2}{\item[({\it ii})]}
\newcommand{\3}{\item[({\it iii})]}
\newcommand{\4}{\item[({\it iv})]}
\newcommand{\5}{\item[({\it v})]}
\newcommand{\6}{\item[({\it vi})]}
\newcommand{\7}{\item[({\it vii})]}
\newcommand{\8}{\item[({\it viii})]}
\newcommand{\9}{\item[({\it xi})]}
\newcommand{\ido}{\int_{\partial\Omega}}
\newcommand{\half}{\mbox{$\frac{1}{2}$}}
\def\parallel{\|}
\def\mid{|}
\def\Bbb R{\real}
\def\hat{\widehat}
\def\tilde{\widetilde}
\def\bar{\overline}
\newcommand{\threehalves}{3\over 2}
\newcommand{\bfPi}{\mbox{\boldmath $\Pi$}}
\newcommand{\bfXi}{\mbox{\boldmath $\Xi$}}
\newcommand{\bfalpha}{\mbox{\boldmath $\alpha$}}
\newcommand{\bfbeta}{\mbox{\boldmath $\beta$}}
\newcommand{\bfgamma}{\mbox{\boldmath $\gamma$}}
\newcommand{\bfdelta}{\mbox{\boldmath $\delta$}}
\newcommand{\bfzeta}{\mbox{\boldmath $\zeta$}}
\newcommand{\bfUpsilon}{\mbox{\boldmath $\Upsilon$}}
\newcommand{\bfGamma}{\mbox{\boldmath $\Gamma$}}
\newcommand{\bfcala}{\mbox{\boldmath ${\cal A}$}}
\newcommand{\bfcalm}{\mbox{\boldmath ${\cal M}$}}
\newcommand{\bfcaln}{\mbox{\boldmath ${\cal N}$}}
\newcommand{\bfcalq}{\mbox{\boldmath ${\cal Q}$}}
\newcommand{\bfcalb}{\mbox{\boldmath ${\cal B}$}}
\newcommand{\bfcalc}{\mbox{\boldmath ${\cal C}$}}
\newcommand{\bfcali}{\mbox{\boldmath ${\cal I}$}}
\newcommand{\bfcalg}{\mbox{\boldmath ${\cal G}$}}
\newcommand{\bfcalh}{\mbox{\boldmath ${\cal H}$}}
\newcommand{\bfcalk}{\mbox{\boldmath ${\cal K}$}}
\newcommand{\bfcalt}{\mbox{\boldmath ${\cal T}$}}
\newcommand{\bfcalx}{\mbox{\boldmath ${\cal X}$}}
\newcommand{\bfcall}{\mbox{\boldmath ${\cal L}$}}
\newcommand{\bfcalf}{\mbox{\boldmath ${\cal F}$}}
\newcommand{\bfcalr}{\mbox{\boldmath ${\cal R}$}}
\newcommand{\bfcals}{\mbox{\boldmath ${\cal S}$}}
\newcommand{\bfcalw}{\mbox{\boldmath ${\cal W}$}}
\newcommand{\bfcalu}{\mbox{\boldmath ${\cal U}$}}
\newcommand{\bfcalv}{\mbox{\boldmath ${\cal V}$}}
\newcommand{\bfcalz}{\mbox{\boldmath ${\cal Z}$}}
\pagenumbering{roman}
\newcommand{\art}[6]{{\I[{\sc #1,}] {#2}, {\it #3}, {\bf #4}, {#5} {[#6]}}}
\newcommand{\ED}{\end{description}}
\newcommand{\I}{\item }
\newcommand{\ra}{\rm a}
\newcommand{\rb}{\rm b}
\newcommand{\rc}{\rm c}
\newcommand{\Hsp}{{\rm I}\!\!\,{\rm R}^n_+}
\newcommand{\Hsn}{{\rm I}\!\!\,{\rm R}^n_-}
\newcommand{\po}[1]{\mbox{$\displaystyle \frac{\mbox{$\partial #1$}}
{\mbox{$\partial x_{1}$}}$}}
\newcommand{\PO}[1]{\mbox{$\displaystyle \frac{\mbox{$\partial #1$}}
{\mbox{$\partial y_{1}$}}$}}
\newcommand{\OP}{\left(\Delta+2\lambda\PO{}\right)}
\newcommand{\op}{\left(\Delta+2\lambda\po{}\right)}
\newcommand{\ft}[1]{
\Frac{1}{(2\pi)^{n/2}}\Int{{\Bbb R}^{n}}{}e^{i{\bf x}\cdot \bfxi}
#1(\xi)d\xi}
\newcommand{\Ft}[1]{
\Frac{1}{2\pi}\Int{{\Bbb R}^{2}}{}e^{i{x}\cdot \xi}
#1(\xi)d\xi}
\newcommand{\Z}{\item[({\it a})]}
\newcommand{\B}{\item[({\it b})]}
\newcommand{\C}{\item[({\it c})]}
\newcommand{\D}{\item[({\it d})]}
\newcommand{\E}{\item[({\it e})]}
\newcommand{\G}{\item[({\it g})]}
\newcommand{\Š}{\`e}
\newcommand{\…}{\`a}
\newcommand{\•}{\`o}
\newcommand{\—}{\`u}
\newcommand{\}{\`{\i}}
\def\tag{\renewcommand{\theequation}}
\newcommand{\Footnote}{~\footnote}
\newcommand{\ie}{{\it i.e.}}
\newcommand{\dist}{\mbox{\rm dist\,}}
\newcommand{\const}{\mbox{\rm const}}
\newcommand{\trace}{\mbox{\rm trace}}
\newcommand{\Bo}{\par\hfill{$\Box$}\par\noindent}
\newcommand{\Nor}[1]{\langle{#1}\rangle_q}
\newcommand{\vs}{\vspace*{.5cm}\par\noindent}
\newcommand{\Vs}{\vspace*{.6cm}\par\noindent}
\newcommand{\Vvs}{\vspace*{.7cm}\par\noindent}
\newcommand{\VVs}{\vspace*{.8cm}\par\noindent}
\newtheorem{definition}{Definition}[section]
\newcommand{\Bd}{\begin{definition}\begin{rm}}
\newcommand{\Ed}{\end{rm}\end{definition}}
\newtheorem{remark}{Remark}[section]
\newcommand{\Br}{\begin{remark}\begin{rm}}
\newcommand{\Er}{\end{rm}\end{remark}}
\newtheorem{proposition}{Proposition}[section]
\newcommand{\Bp}{\begin{proposition}\begin{sl}}
\newcommand{\EP}[1]{\end{sl}\label{proposition:#1}\end{proposition}}
\newcommand{\propref}[1]{{\rm Proposition \ref{proposition:#1}}}
\newcommand{\Bt}{\begin{theorem}\begin{sl}}
\newcommand{\Et}{\end{sl}\end{theorem}}
\newcommand{\Bl}{\begin{lemma}\begin{sl}}
\newcommand{\El}{\end{sl}\end{lemma}}
\newtheorem{theorem}{Theorem}[section]
\newtheorem{lemma}{Lemma}[section]
\newtheorem{corollary}{Corollary}[section]
\newcommand{\eqref}[1]{{\rm (\ref{eq:#1})}}
\newcommand{\Bc}{\begin{corollary}\begin{sl}}
\newcommand{\Ec}{\end{sl}\end{corollary}}
\newcommand{\ET}[1]{\end{sl}\label{theorem:#1}\end{theorem}}
\newcommand{\EDD}[1]{\end{rm}\label{definition:#1}\end{definition}}
\newcommand{\EL}[1]{\end{sl}\label{lemma:#1}\end{lemma}}
\newcommand{\theoref}[1]{{\rm Theorem \ref{theorem:#1}}}
\newcommand{\defref}[1]{{\rm Definition \ref{definition:#1}}}
\newcommand{\ER}[1]{\end{rm}\label{remark:#1}\end{remark}}
\newcommand{\EC}[1]{\end{sl}\label{corollary:#1}\end{corollary}}
\newcommand{\remref}[1]{{\rm Remark \ref{remark:#1}}}
\newcommand{\cororef}[1]{{\rm Corollary \ref{corollary:#1}}}
\newcommand{\lemmref}[1]{{\rm Lemma \ref{lemma:#1}}}
\newcommand{\essup}[1]{{\rm ess}\,{{\displaystyle \sup_{\hspace*{-5mm}{#1}}}}}

\renewcommand{\i}{{\rm i}}

\pagenumbering{arabic}
\newcommand{\QED}{{\par\hfill$\square$\par}}
\renewcommand{\thefootnote}{(\arabic{footnote})}
\title{Very Weak Solutions and  Asymptotic Behavior of Leray Solutions\\ to the Stationary  Navier-Stokes Equations}
\author{ Giovanni P. Galdi 
\thanks{Department of Mechanical Engineering and Materials Science, University of Pittsburgh, PA 15261. 
}}
\date{}
\maketitle
\begin{abstract}  Let $\bfu$ be a Leray solution to the Navier-Stokes boundary-value problem in an exterior domain, vanishing at infinity and satisfying the generalized energy inequality. We show that if there exist $R>0$ and ${\sf s}\ge \frac23 q$, $q>6$, such that the $L^{\sf s}-$norm of $\bfu$ on the spherical surface of radius $R$ divided by $R$ is less than a constant depending only on {\sf s} and $q$, then $\bfu(x)$ must decay as $|x|^{-1}$ for $|x|\to\infty$. This result is proved with an approach based on a new theory of very weak solutions in exterior domains which, as such, is of independent interest.     
 \end{abstract}

\renewcommand{\theequation}{\arabic{section}.\arabic{equation}}
\setcounter{section}{0}
\section{Introduction} One of the still unresolved issues in the theory of the Navier-Stokes boundary-value problem in a domain, $\Omega$, exterior to a smooth compact set of $\real^3$ (flow around a body), concerns the behavior at large spatial distances of Leray solutions. We recall that the latter are characterized by a velocity field $\bfu\in L^6(\Omega)$ with $\nabla\bfu\in L^2(\Omega)$  satisfying (in an appropriate sense) the following set of equations
\be    
\ba{cc}\medskip\left.\ba{ll}\medskip
\Delta\bfu=\nabla {\sf p}+\bfu\cdot\nabla\bfu+\bfF\\ \medskip
\Div\bfu=0\ea\right\}\ \ \mbox{in $\Omega$\,,}\\
\bfu=\bfu_*\ \ \mbox{at $\partial\Omega$}\,,\ \ \Lim{|x|\to\infty}\bfu(x)=\0\,,
\ea
\eeq{1.1}
where $\bfF$ and $\bfu_*$ are prescribed functions. 
A classical result of Leray \cite{Leray1,Leray2} (see also \cite{Finn} ), ensures  that the class of such  solutions is not empty, provided only $\bfF$, $\bfu_*$ and $\partial\Omega$  are sufficiently regular and the magnitude of the flux of $\bfu_*$ through $\partial\Omega$ is suitably restricted. Moreover, it is possible to associate to $\bfu$ a scalar field ${\sf p}$ such that, for smooth enough data, $(\bfu,{\sf p})$ satisfies \eqref{1.1} in the ordinary sense.
It is also a well-established fact \cite[Theorem IX.6.1]{Gab0} that, if $\bfF$ has a bounded support or, more generally,  $\bfF$ and all its derivatives, for large $|x|$, belong to appropriate Lebesgue spaces, then
\be
\lim_{|x|\to\infty}D^\alpha\bfu(x)=\0\,,\ \ \mbox{all $|\alpha|\ge0$\,.}
\eeq{1.3}
\par
The question that has intrigued mathematicians since Leray's pioneering work and which still remains unanswered is whether the conditions \eqref{1.3} are satisfied with a {\em specific} and {\em sharp} decay rate. Investigating this question is not just a matter of mathematical completeness, as it is closely related to the validity of fundamental physical properties, which any solution representing a fluid flow around a body should possess; see \cite[Introduction to X]{Gab}. Moreover, the knowledge of the sharp asymptotic behavior may play a fundamental role in the resolution of the Liouville problem \cite[p. 12]{Gab} that, in turn, is tightly related to existence of possible singularities for the Cauchy problem \cite{KNSS,AB}. 
It should also be noted that the problem is unresolved only if $\bfu$ vanishes for $|x|\to\infty$. Indeed, if $\bfu$ converges to a non-zero vector $\bfu_\infty$, or even to the velocity field  $\bfu_\infty+\bfomega\times\bfx$, $\bfomega\in\real^3$, $\bfu_\infty\neq\0$, of a generic rigid motion, one is able to provide a complete and distinctive asymptotic behavior; see \cite[Chapters X and XI]{Gab}, \cite{Kyed}.
\par
The first result devoted to the problem was furnished by the author  \cite{GaDD} when $\bfu_*\equiv\0$, in the class of Leray solutions satisfying the {\em energy inequality}:
\be
\int_{\Omega}|\nabla\bfu|^2\le \int_\Omega\bfF\cdot\bfu\,.
\eeq{1.2}
As is well known, this class is not empty (e.g., \cite[Theorem X.4.1]{Gab}). Specifically, in \cite{GaDD} it was shown that 
if the magnitude of $\bfF$ in the $L^1$-norm and in a suitably weighted $L^\infty$-norm  is below a certain constant, then the corresponding Leray solution obeying \eqref{1.2} 
 decays like $|x|^{-1}$ as $|x|\to\infty$. It was later proved \cite{KSv}  that, in fact, if $\bfF$ is of bounded support, $\bfu$ behaves asymptotically precisely like a suitable Landau solution, which implies that  $O(|x|^{-1})$ is sharp; see also \cite[Remark X.9.3]{Gab}. The result in \cite{GaDD} was generalized to higher dimensions in \cite{Miy} (see also \cite{KS}). More recently, further remarkable contributions to the question were furnished in \cite{Car1,Car2,ChBj,Koz,Weng1,Weng2} under the assumptions of axial-symmetry. There,  the authors are able to show  suitable algebraic decay rate for $\bfu$ and some of its derivatives,  in the direction orthogonal to the axis of symmetry. It must be noticed that in these works, the validity of the energy inequality is not needed.
\par
The aim of this paper is to provide an additional contribution to the problem.
Precisely, let $\bfu$ be a Leray solution to \eqref{1.1} with $\bfF$ of bounded support, satisfying the generalized energy inequality (see \eqref{4.2}).\footnote{If $\bfu_*\equiv\0$ this inequality reduces to \eqref{1.3}.} With respect to
 spherical coordinates
$(|x|,\omega)$  with the origin in the interior of $\real^3\backslash\Omega$, denote by $\bfu(R,\omega)$ the trace of $\bfu$ on the sphere $|x|=R$. Moreover, let $q>6$ and ${\sf s}\ge \frac23 q$. We then show that there exists a constant $C_0=C_0({\sf s},q)$ such that if ($S^2\equiv$ unit sphere)
\be
R\,\left(\int_{S^2}|\bfu(R,\omega)|^{\sf s}{\rm d}\sigma_\omega\right)^\frac1{\sf s}\le C_0\,,\ \ \,\mbox{for {\em some} sufficiently large $R$}\,,  
\eeq{1.4}
necessarily $\bfu(x)=O(|x|^{-1})$ as $|x|\to\infty$. From known results, this then implies also $D^\alpha\bfu(x)=O(|x|^{-1-|\alpha|})$, for all $|\alpha|\ge 0$.
\par  
The proof of this statement is obtained by means of an  approach based on the theory of very weak solutions. As is known, these solutions are characterized by being generated from boundary and volume data of very low regularity, only belonging to appropriate negative Sobolev spaces. They are ``very weak" because, at the outset, they do not possess any locally
integrable derivative and satisfy the relevant equations only in a  distributional sense. In the stationary context, they were introduced in \cite{Giga} for the linear (Stokes) problem and were further analyzed and completed in \cite{GSS}, including the full nonlinear case (see also \cite{Conca,MP}). All these papers deal with solutions in bounded domains. The extension to the case at hand of an exterior domain was given in \cite{Kim} (see also \cite{CoTaRu}) but, however, under the crucial assumption that $\bfu(x)$ converges, for $|x|\to\infty$, to a given {\em non-zero} vector, which is not the case considered here. Therefore, our first objective is to  develop a full theory of existence and uniqueness of very weak solutions to \eqref{1.1}. This is done by studying first  this theory for the linear (Stokes) problem, obtained by disregarding the nonlinearity in \eqref{1.1}$_1$; see Theorems 2.1 and 2.2. Successively, by means of an appropriate perturbation argument, results are extended to the full nonlinear case in \theoref{3}. In particular, under the assumption that $\bfF:=\Div\mathbb H$, this theorem ensures  that if  $\bfu_*\in W^{-\frac1q,q}(\partial\Omega)$, $q>6$, $\mathbb H\in L^{\frac{q}2}$ ``near" the boundary and decays, for large $|x|$, like $|x|^{-2}$, then for data of restricted magnitude there is a unique corresponding very weak solution to \eqref{1.1} that, for large $|x|$, decays like $|x|^{-1}$. 
\par
The basic idea that connects the study of the asymptotic behavior of Leray solutions with the theory of very weak solutions is the following.
Let $\Omega^R:=\{x\in\Omega:\,|x|>R\}$ where $R$ is so large that the  support of $\bfF$ is contained in the interior of the complement of $\Omega^R$. Let $\bfu(R,\omega)=:\bfv_*(\omega)$ be the trace of a Leray solution $\bfu$ on $\partial\Omega_R$. From classical theorems on the regularity of Leray solutions we have that $\bfv_*\in C^{\infty}(\partial\Omega_R)$. Thus, clearly, 
\be
\|\bfv_*\|_{{\sf s},{S^2}}:=\left(\int_{S^2}|\bfv_*|^{\sf s}{\rm d}\sigma_\omega\right)^\frac1{\sf s}<\infty\,.
\eeq{1.5}
Imposing {\em only} this regularity condition we then show,  thanks to the results proved in \theoref{3}, that if, in addition, $\|\bfv_*\|_{{\sf s},S^2}$ satisfies a restriction as in \eqref{1.4},  
there exists a corresponding very weak solution, $\bfv=\bfv(x)$, to \eqref{1.1} with $\bfF\equiv\0$, $\bfu_*\equiv\bfv_*$ and $\Omega\equiv\Omega^{R}$. However, since $\bfv_*$ and $\Omega^R$ are both $C^\infty$, it turns out that $\bfv$ is $C^\infty$ as well, up to the boundary. Consequently, also using the property that $\bfv$ decays as $|x|^{-1}$, we show that it further obeys the energy equality; see \theoref{3.2}. Hence, since $\bfu$ satisfies the generalized energy inequality, we can use a  weak-strong uniqueness argument similar to that employed in \cite{GaDD} to secure  that $\bfv$ must coincide with  $\bfu$ in $\Omega^{R}$, thus recovering the desired result; see \theoref{4.1}.
Although the assumption \eqref{1.4} appears rather mild, it remains an open question whether it is satisfied by Leray solutions that obey the generalized energy inequality. Indeed, to date, we only know that such solutions satisfy \eqref{1.4} with $R$ replaced by $R^{\frac12}$; see \remref{4.1}.  
\par
The paper outline is as follows. In Section 2 we begin the study of very weak solutions to the Stokes problem, obtained from \eqref{1.1} by omitting in it the nonlinear term $\bfu\cdot\nabla\bfu$. Precisely, after giving the definition of very weak solution, we show their uniqueness in \theoref{0} and existence in \theoref{1}. The important properties of these solutions is that they become regular away from the boundary and, in particular, decay as $|x|^{-1}$ for large $|x|$. Combining these results with a contraction map argument, in Section 3 we extend the existence and uniqueness theory to the full nonlinear case, provided the size of the data is suitably restricted; see \theoref{3}. By using a scaling argument, we then specialize  the finding of \theoref{3} to the case where $\Omega$ is the exterior of a ball of radius $R$ and $\bfF\equiv\0$, and show, in particular, existence under the sole condition that the boundary data are in $L^{\sf s}(S^2)$ and satisfy a condition of the type \eqref{1.4}; see \theoref{3.2}. 
In the final Section 4, we give the precise definition of Leray solution $\bfu$ satisfying the generalized energy inequality and recall in \propref{4.1} a corresponding existence result. Successively, we combine the result of \theoref{3.2} with the weak-strong uniqueness argument used in \cite{GaDD} to prove the desired decay property for $\bfu$ in \theoref{4.1}. The paper ends with an Appendix dedicated to the existence and estimate of solutions to a Stokes problem in the whole of $\real^3$, a technical result that is needed in the proof of \theoref{1}.

\setcounter{equation}{0}
\section{Existence and uniqueness of very weak solutions to the Stokes problem }
We begin by recalling some notation. An exterior domain, $\mathcal A$, is a domain of $\real^3$, whose complement $\mathcal A^c$ coincides with the closure of a bounded domain. We set $\delta_\cala={\rm diam}\,\cala^c$.  By $\Omega$  we will indicate an exterior  domain  of class $C^2$, even though such a regularity assumption is not always necessary. For simplicity, we put $\delta\equiv\delta_{\Omega}$.\footnote{We remark that all our results continue to hold if $\Omega\equiv\real^3$.} 
We take the origin of coordinates in the interior of $\cala^c$, and denote by $(|x|,\omega:=(\theta,\varphi))$ the corresponding spherical coordinates, and ${\rm d}\,\sigma_\omega$ the infinitesimal solid angle. Let  $B_r:=\{x\in\real^3:\,|x|<r\,,\ r>0\}$, $S^2:=\partial B_1$, and set 
$$\cala_R:=\cala\cap B_R,\  R>\delta_{\cala}\,,\ \
\cala^R:=\cala\backslash\bar{\cala_R}\,,\ \  \cala_{R_1,R_2}:=\cala^{R_1}\backslash\bar{\cala^{R_2}}\,, \ \ R,R_1,R_2>\delta_\cala\,, \ \ R_2>R_1\,.
$$ 
For $A$ a domain of $\real^3$,  $L^q=L^q(A)$,  $q\in[1,\infty]$,  denote  Lebesgue spaces,  $\|\cdot\|_{q,A}$ their norm and $\langle\cdot ,\cdot \rangle_A$  the $L^{q'}(A)-L^q(A)$ duality pairing, $\frac1{q'}+\frac1q=1$ (scalar product if $q=2$). Moreover,
$W^{m,q}=W^{m,q}(A)$, $m\in\nat$, are  Sobolev spaces with norm  $\|\cdot\|_{m,q,A}$, and $W^{-m,q'}(A)$  their dual, with norm $\|\cdot\|_{-m,q^\prime,A}$. Trace spaces are denoted by $W^{m-\frac1q,\frac1q}(\partial A)$, with  norm  $\|\cdot\|_{m-\frac1q,\frac1q,\partial A}$, and associated dual $W^{-m+\frac1{q'},\frac1{q'}}(\partial A)$, with norm $\|\cdot\|_{-m+\frac1{q'},\frac1{q'},\partial A}$. The duality pairing between trace spaces will be indicated by $\langle\cdot,\cdot\rangle_{\partial A}$.  Moreover, $D^{m,q}=D^{m,q}(A)$ denotes the homogeneous Sobolev space with semi-norm $\sum_{|l|=m}\|D^lu\|_{q,A}$, while $D_0^{m,q}(A)$ is the completion of $C_0^\infty(A)$ in the norm $\sum_{|l|=m}\|D^lu\|_{q,A}$. The dual space of $D_0^{m,q}(A)$ will be indicated with $D_0^{-m,q'}(A)$. In all the above notation, the subscript ${}_A$ will be omitted, unless confusion may arise.  
\smallskip\par
Consider the following Stokes boundary-value problem
\be\ba{cc}\medskip\left.\ba{ll}\medskip
\Delta\bfv=\nabla p+\Div\mathbb F\\ \medskip
\Div\bfv=0\ea\right\}\ \ \mbox{in $\Omega$\,,}\\
\bfv=\bfv_*\ \ \mbox{at $\partial\Omega$}\,,\ \ \Lim{|x|\to\infty}\bfv(x)=\0\,,
\ea
\eeq{1}
where the second order tensor  $\mathbb F$ and the vector  $\bfv_*$ are assigned fields. 
The main objective of this section is to prove the existence,  uniqueness and associated estimates of solutions to \eqref{1} when $\bfv_*$ is in spaces of low regularity (negative Sobolev spaces). 
\par
To this end, denote by $\calk_0^2=\calk_0^2(\Omega)$ the class of functions $\bfphi\in C^2(\Omega)$ such that
\begin{itemize}
  \item [(i)]
$\bfphi\in C^2(\bar{\Omega_R})\ \mbox{for all $R>\delta$}$\,;
\item[(ii)] $\bfphi=0$\ \mbox{at $\partial\Omega$}\,;
\item[(iii)] $\bfphi(x)=0$\ \mbox{for all $|x|\ge \rho_{\varphi}$, some $\rho_{\varphi}>\delta$}.  
\end{itemize}
Following \cite{GSS}, we now give the definition of very weak solution to \eqref{1}.
\Bd
The field $\bfv:\Omega\to \real^3$ is a very weak solution to \eqref{1} if for some $q\in (1,\infty)$ the following conditions are fulfilled.
\begin{itemize}
  \item [(a)] $\bfv\in L^q(\Omega_R)$, for all $R>\delta$\,;
  \item [(b)] $\bfv$ satisfies
 $$ 
\langle\bfv,\Delta\bfphi\rangle+\langle\bfv_*,\bfn\cdot{\nabla\bfphi}\rangle_{\partial\Omega}+\langle\mathbb F,\nabla\bfphi\rangle=0\,,\ \,\mbox{for all $\bfphi\in \calk_0^2(\Omega)$}\,;
$$
\item[(c)]\ $\langle\bfv,\nabla\zeta\rangle=0$\,,\ \ \mbox{for all $\zeta\in C_0^1(\real^3)$}\,;
\item[(d)]\ $\Liminf{|x|\to\infty}\int_{S^2}|\bfv(|x|,\omega)|{\rm d}\sigma_\omega=0\,.$ 
\end{itemize}
\EDD{1}
In the next two theorems we will furnish uniqueness and existence of very weak solutions under suitable assumptions on the data. We begin with a uniqueness result.\footnote{The assumptions of \theoref{0} can be significantly weakened. However, this would be irrelevant for our purposes.}
\Bt Let $\bfv$ be a very weak solution corresponding to $\mathbb F\equiv\bfv_*\equiv\0$. Then, if for some ${\sf R}>0$ and $q>1$, 
$$
\bfv\in D^{1,q}(\Omega^{\sf R})\,,
$$
necessarily $\bfv\equiv\0$.
\ET{0}
{\em Proof.} By assumption and  definition of very weak solution,  we infer that $\bfv\in W^{1,q}(\Omega_{{\sf R},2{\sf R}})$ which, by known trace theorems, entails in particular that $\bfv\in W^{1-\frac1q,q}(\partial B_{2{\sf R}})$. Consequently, $\bfv$ is a very weak solution in the bounded domain $\Omega_{2{\sf R}}$, with boundary data in $W^{1-\frac1q,q}(\partial \Omega_{2{\sf R}})$ and $\mathbb F\equiv\0$. So,    
from \cite[Lemma 4]{GSS}, it follows that $\bfv\in W^{1,q}(\Omega_{2{\sf R}})$. Combining the latter with the assumption, we then obtain that $\bfv$ is a $q$-generalized solution to \eqref{1}; see \cite[Definition V.1.1]{Gab}. Therefore, by \cite[Theorem V.3.4]{Gab}, it must vanish identically in $\Omega$.
\par\hfill$\square$
\par
 The following existence result holds.
\Bt Let $q\in(\frac32,\infty)$, $\rho>s>\delta$. Assume that $\mathbb F$ and $\bfv_*$ satisfy
$$
\big(|x|^2\mathbb F\big)\in L^\infty(\Omega^s)\,;\ \ \mathbb F\in L^r(\Omega_{s})\,,\ \ \mbox{$\frac1r\le \frac1q+\frac13\,,\ r\in (1,q]$}\,;\ \ \bfv_*\in W^{-\frac1q,q}(\partial\Omega).
$$
Then, there is  a unique corresponding very weak solution such that
$$
(|x|\,\bfv)\in L^\infty(\Omega^{2\rho})\,,\ \bfv\in D^{1,q}(\Omega^\rho)\cap W^{1,r}(\Omega_{s,\rho})\,.
$$ 
Furthermore, there exists $p\in L^q(\Omega^\rho)\cap W^{-1,q}(\Omega_{\rho})$ such that $(\bfv,p)$ satisfies \eqref{1}$_1$ in the sense of distributions.
Finally, the  following estimate holds
\be\ba{rl}\medskip
\|\big(|x|\,\bfv\big)\|_{\infty,\Omega^{2\rho}}+\|\nabla\bfv\|_{q,\Omega^\rho}+\|\bfv\|_{1,r,\Omega_{s,\rho}}+&\!\!\!\!\|\bfv\|_{q,\Omega_{2\rho}}+\|p\|_{q,\Omega^\rho}+\|p\|_{-1,q,\Omega_{\rho}}\\
&\!\!\!\!\le c\left(\|\big(|x|^2\,\mathbb F\big)\|_{\infty,\Omega^s}+\|\mathbb F\|_{r,\Omega_{s}}+\|\bfv_*\|_{-\frac1q,q,\partial\Omega}\right)\,,
\ea
\eeq{2}
where $c=c(\Omega,\rho,s,r,q)>0$. 
\ET{1}
{\em Proof.} Uniqueness is a direct consequence of \theoref{0}. To prove existence,   
assume first $\bfv_*\in W^{1-\frac1q,q}(\partial\Omega)$. We extend $\mathbb F$ to 0 in $\Omega^c$ (continuing to denote the extension by $\mathbb F$) and for each $\varepsilon:=\frac1m\in(0,1]$, $m\in\nat$, we define the (Friedrich) mollifier of $\mathbb F$:
$$
\mathbb F_\varepsilon(x):=\int_{\real^3}k_\varepsilon(x-y)\mathbb F(y){\rm d}y\,,\ \ k_\varepsilon(\xi):=\varepsilon^{-3}k(\xi/\varepsilon)\,,\ \ k\in C_0^\infty(B_1)\,,
 \ \ \int_{\real^3}k(\xi){\rm d}\xi=1.
$$
Setting 
\be
\hat{\mathbb F}_\varepsilon(x):=\left\{\ba{ll}\medskip 
\mathbb F(x)\ \mbox{if $x\in\Omega^s$}\\
\mathbb F_\varepsilon(x)\ \mbox{if $x\in\Omega_s$}
\ea\right.\,,
\eeq{0}
from the assumption on $\mathbb F$ and elementary properties of mollifiers, we deduce
$$
\big(|x|^2\hat{\mathbb F}_\varepsilon\big)\in L^\infty(\Omega)\,.
$$
Then, by \cite[Theorem V.8.1 and Exercise V.8.1]{Gab}, there exists a unique $q$-generalized solution $(\bfv_\varepsilon,p_\varepsilon)$ to \eqref{1} with $\mathbb F \equiv \hat{\mathbb F}_\varepsilon$, such that
$$
(|x|\,\bfv_\varepsilon)\in L^\infty(\Omega^R)\,,\ (\bfv_\varepsilon,p_\varepsilon)\in D^{1,q}(\Omega)\times L^q(\Omega)\,,
$$ 
where $R>\delta$ is arbitrary.
Let $\psi=\psi(|x|)$ be a non-decreasing smooth function such that $\psi(|x|)=0$, for $|x|\in [0,s]$, and $\psi(|x|)=1$ for $x\in \Omega^\rho$, and set
\be
\bfu_\varepsilon:=\psi\,\bfv_\varepsilon\,,\ \ {\sf p}_\varepsilon:=\psi\,p_\varepsilon\,,
\eeq{3} 
from \eqref{1},  \eqref{0} and the definition of $\psi$, we deduce that $(\bfu_\varepsilon,{\sf p}_\varepsilon)$ solves the following problem
\be\medskip\left.\ba{ll}\medskip
\Delta\bfu_\varepsilon=\nabla {\sf p}_\varepsilon+\Div\mathbb G+\bff_\varepsilon\\ \medskip
\Div\bfu_\varepsilon=g_\varepsilon\ea\right\}\ \ \mbox{in $\real^3$\,,}
\eeq{4}
where
\be\ba{rl}\medskip
(\mathbb G)_{ij}:=&\!\!\!\psi\, (\mathbb F)_{ij}\\
\medskip
f_{\varepsilon i}:=&\!\!\![(\mathbb D(\bfv_\varepsilon))_{ik}-{p_\varepsilon}\,\delta_{ik}]D_k\psi+D_k(v_{\varepsilon k}D_i\psi+
v_{\varepsilon i}D_k\psi)-(\hat{\mathbb F}_{\varepsilon})_{ ki}D_k\psi\,,\\ \medskip
\mathbb D(\bfv_\varepsilon):=&\!\!\!\half(\nabla\bfv_\varepsilon+(\nabla\bfv_\varepsilon)^\top)\,,
\\
g_\varepsilon:=&\!\!\!\bfv_\varepsilon\cdot\nabla\psi\,.
\ea\eeq{5}
The solution $(\bfu_\varepsilon,{\sf p}_\varepsilon)$ obeys the following estimate
\be
\|(|x|\,\bfu_\varepsilon)\|_{\infty,B^{2\rho}}+\|\nabla\bfu_\varepsilon\|_{q,\real^3}+\|{\sf p}_\varepsilon\|
_{q,\real^3}
\le c\left(\|(|x|^2\mathbb G)\|_{\infty,\real^3}+\|\bff_\varepsilon\|_{-1,q, B_{2\rho}}+
\|g_\varepsilon\|_{q,B_{\rho}}\right)\,.
\eeq{6}
The validity of this property can be deduced from the proof used for \cite[Lemma V.8.2]{Gab}. However, for the reader's sake, we include a full proof in the Appendix.  
By the assumptions on $q$ and $r$, we have, on the one hand, $1-\frac1q=:\frac1{q'}>\frac13$, and hence, on the other hand, the embedding $W^{1,q'}(B_{2\rho})\subset L^{r'}(B_{2\rho})$, $\frac1{r'}=1-\frac1r$. As a result,  we readily deduce
$$
\|\bff_\varepsilon\|_{-1,q, B_{2\rho}}:=\sup_{\mbox{\footnotesize $\bfphi\in W_0^{1,q'}(B_{2\rho});\,\|\bfphi\|_{1,q'}=1$}}|\langle \bff_\varepsilon,\bfphi\rangle|\le c\left(\|\hat{\mathbb F}_\varepsilon\|_{r,\Omega_{2\rho}}+\|\bfv_\varepsilon\|_{q,\Omega_{2\rho}}+\|p_\varepsilon\|_{-1,q,\Omega_{2\rho}}\right)\,.
$$
Moreover,
$$
\|g_\varepsilon\|_{q,B_{\rho}}\le c\,\|\bfv_\varepsilon\|_{q,\Omega_{\rho}}\,.
$$
Thus, in view of the latter two displayed relations,  \eqref{3}  and \eqref{6} we infer
\be\ba{rl}\medskip
\|\big(|x|\,\bfv_\varepsilon\big)\|_{\infty,\Omega^{2\rho}}+\|\nabla\bfv_\varepsilon\|_{q,\Omega^\rho}+&\!\!\!\!\|{p}_\varepsilon\|
_{q,\Omega^\rho}\\
&\!\!\!\!\le c\left(\|\big(|x|^2\mathbb F\big)\|_{\infty,\Omega^s}+\|\hat{\mathbb F}_\varepsilon\|_{r,\Omega_{2\rho}}+
\|\bfv_\varepsilon\|_{q,\Omega_{2\rho}}+\|p_\varepsilon\|_{-1,q,\Omega_{2\rho}}\right)\,.\ea
\eeq{7}
Next, we observe that, since $\bfv_\varepsilon\in L^q(\Omega_R)$ for all $R>\delta$, by \cite[Theorem 2]{GSS} $\bfv_\varepsilon$ is the unique very weak solution to the following Stokes problem in the bounded domain $\Omega_{2\rho}$:
\be\ba{cc}\medskip\left.\ba{ll}\medskip
\Delta\bfv_\varepsilon=\nabla p_\varepsilon+\Div\hat{\mathbb F}_\varepsilon\\ \medskip
\Div\bfv_\varepsilon=0\ea\right\}\ \ \mbox{in $\Omega_{2\rho}$\,,}\\
\bfv_\varepsilon=\bfv_*\ \ \mbox{at $\partial\Omega$}\,,\ \ \bfv_\varepsilon=\bfv_\varepsilon\ \ \mbox{at $\partial B_{2\rho}$}\,.
\ea
\eeq{8}
As a result, again by  \cite[Theorem 2]{GSS}, it follows that $(\bfv_\varepsilon,p_\varepsilon)$ obeys the estimate
\be
\|\bfv_\varepsilon\|_{q,\Omega_{2\rho}}+\|p_\varepsilon\|_{-1,q,\Omega_{2\rho}}\le c\left(\|\hat{\mathbb F}_\varepsilon\|_{r,\Omega_{2\rho}}+\|\bfv_*\|_{-1,q,\partial\Omega}+\|\bfv_\varepsilon\|_{-1,q,\partial B_{2\rho}}\right)\,.
\eeq{9}
Combining \eqref{7} and \eqref{9} we then conclude, in particular,
\be\ba{rl}\medskip
\|\big(|x|\,\bfv_\varepsilon\big)\|_{\infty,\Omega^{2\rho}}+\|\nabla\bfv_\varepsilon\|_{q,\Omega^\rho}&\!\!\!\!+\|\bfv_\varepsilon\|_{q,\Omega_{2\rho}}+\|p_\varepsilon\|_{q,\Omega^\rho}+\|p_\varepsilon\|_{-1,q,\Omega_\rho}\\
&\!\!\!\le c\left(\|\big(|x|^2\mathbb F\big)\|_{\infty,\Omega^s}+\|\hat{\mathbb F}_\varepsilon\|_{r,\Omega_{2\rho}}+\|\bfv_*\|_{-\frac1q,q,\partial\Omega}+\|\bfv_\varepsilon\|_{-\frac1q,q,\partial B_{2\rho}}\right)\,,
\ea
\eeq{10}
where $c=c(\Omega,\rho,r,q)>0$.   We want to show that there exists $c=c(\Omega,\rho,r,q)>0$ such that
\be
\|\bfv_\varepsilon\|_{-\frac1q,q,\partial B_{2\rho}}\le c\left(\|\psi\,\big(|x|^2\mathbb F\big)\|_{\infty,\Omega}+\|\hat{\mathbb F}_\varepsilon\|_{r,\Omega_{2\rho}}+\|\bfv_*\|_{-\frac1q,q,\partial\Omega}\right)\,,
\eeq{11}
for all $\varepsilon\equiv \frac1m\in (0,1]$, $m\in\nat$.
The proof relies on a classical contradiction argument based on uniqueness and local compactness. Thus, suppose that there exist sequences of data $\{\mathbb F_n\}$ and $\{\bfv_{*n}\}$ such that, denoting by $\{\bfv_n\}$ the sequence of corresponding very weak solutions, it is 
\be
\|\bfv_n\|_{-\frac1q,q,\partial B_{2\rho}}=1\,,\ \ \|\psi\,\big(|x|^2\mathbb F_n\big)\|_{\infty,\Omega}+\|\mathbb F_n\|_{r,\Omega_{2\rho}}+\|\bfv_{*n}\|_{-\frac1q,q,\partial\Omega}< n^{-1}\,.
\eeq{12}
Because of \eqref{10} and \eqref{12}, we deduce for arbitrary $t>3$
$$
\|\bfv_n\|_{t,\Omega^{2\rho}}+\|\nabla\bfv_n\|_{q,\Omega^\rho}+\|\bfv_n\|_{q,\Omega_{2\rho}}+\|p_n\|_{q,\Omega^\rho}+\|p_n\|_{-1,q,\Omega_\rho}\le C\,,
$$
where the positive constant $C$ is independent of $n$. From the latter, it follows the existence of a pair $(\hat{\bfv},\hat{p})$ with 
$$
(\hat{\bfv},\hat{p})\in [D^{1,q}(\Omega^\rho)\cap L^t(\Omega^\rho)\cap L^q(\Omega_\rho)]\times [L^q(\Omega^\rho)\cap W^{1-\frac1q}(\Omega_\rho)]\,,
$$
and such that (possibly, along a subsequence)
\be\ba{ll}\medskip
\bfv_n\to\hat{\bfv} \  \mbox{weakly in $L^t(\Omega^{2\rho})\cap L^q(\Omega_{2\rho})$}\,,\ \ \nabla\bfv_n\to\nabla\hat{\bfv} \ \ \mbox{weakly in $L^q(\Omega^\rho)$}\,;\\
p_n\to\hat{p} \ \, \mbox{weakly in $L^q(\Omega^\rho)\cap W^{1-\frac1q,q}(\Omega_\rho)$}\,.
\ea
\eeq{13}
It is easy to prove that, in view of \eqref{12}$_2$--\eqref{13}, the field $\hat{\bfv}$ is a very weak solution to \eqref{1} corresponding to $\mathbb F\equiv\bfv_*\equiv\0$. Since $\hat{\bfv}\in D^{1,q}(\Omega^\rho)$, from \theoref{0} we deduce $\hat{\bfv}\equiv\0$. As a result, by \eqref{13}$_{1,2}$  
and a classic compact embedding theorem we obtain, in particular,
\be
\bfv_n\to \0\ \, \mbox{strongly in $L^q(\Omega_{2\rho,3\rho})$}\,.
\eeq{14}
Now, employing a well-known trace inequality \cite[Exercise II.4.1]{Gab} and \eqref{13}$_2$, we obtain, for any $\eta>0$, 
$$
\|\bfv_n\|_{-\frac1q,q,\partial B_{2\rho}}\le \|\bfv_n\|_{q,\partial B_{2\rho}}\le c_\eta\|\bfv_n\|_{q,\Omega_{2\rho,3\rho}}+\eta\,\|\nabla\bfv_n\|_{q,\Omega_{2\rho,3\rho}}\le c_\eta\, \|\bfv_n\|_{q,\Omega_{2\rho,3\rho}}+\eta \,M
$$
with $M$ independent of $n$ and $\eta$. Therefore, letting $n\to\infty$ in the latter and 
taking into account \eqref{14}, we conclude 
$$
\lim_{n\to\infty}\|\bfv_n\|_{-\frac1q,q,\partial B_{2\rho}}=0\,, 
$$
which contradicts \eqref{12}$_1$. Therefore, \eqref{11} is established and so, from this,  \eqref{10}, and the inequality 
\be
\|\hat{\mathbb F}_\varepsilon\|_{r,\Omega_{2\rho}}\le \|\mathbb F\|_{r,\Omega_{2\rho}}\le \|\mathbb F\|_{r,\Omega_{s}}+c\,\|(|x|^2\mathbb F)\|_{\infty,\Omega^s}\,,
\eeq{16} 
we conclude 
\be\ba{rl}\medskip
\|\big(|x|\,\bfv_\varepsilon\big)\|_{\infty,\Omega^{2\rho}}+\|\nabla\bfv_\varepsilon\|_{q,\Omega^\rho}&\!\!\!\!+\|\bfv_\varepsilon\|_{q,\Omega_{2\rho}}+\|p_\varepsilon\|_{q,\Omega^\rho}+\|p_\varepsilon\|_{-1,q,\Omega_\rho}\\
&\!\!\!\le c\left(\|\big(|x|^2\mathbb F\big)\|_{\infty,\Omega^s}+\|\mathbb F\|_{r,\Omega_{s}}+\|\bfv_*\|_{-\frac1q,q,\partial\Omega}\right)\,.
\ea
\eeq{17}
Next, let $\rho'>\rho$. From local estimates on the Stokes problem \cite[Remark IV.4.2]{Gab}, we have 
$$
\| \bfv_\varepsilon\|_{1,r,\Omega_{s,\rho}}\le c\,\big(\|\hat{\mathbb F}_\varepsilon\|_{r,\Omega_{\rho'}}+\|\bfv_\varepsilon\|_{r,\Omega_{\rho'}}\big)\,,
$$
which, in combination with \eqref{16}--\eqref{17} and the assumption $r\le q$, entails
\be
\|\bfv_\varepsilon\|_{1,r,\Omega_{s,\rho}}\le c\,\big(\|{\mathbb F}\|_{r,\Omega_{s}}+\|(|x|^2\mathbb F)\|_{\infty,\Omega^s}\big)\,.
\eeq{18}
Now, clearly, for any given $\varepsilon\equiv \frac1m\in(0,1]$, $m\in\nat$, we have that $\bfv_\varepsilon$ is a very weak solution to \eqref{1} with $\mathbb F\equiv \hat{\mathbb F}_\varepsilon$. Consequently, by passing to the limit $m\to\infty$ and using \eqref{17} and \eqref{18} together with the elementary properties of the mollifiers, one can easily show that the sequence $\{\bfv_\varepsilon\}$ converges to a very weak solution that satisfies all the properties stated in the theorem. Similarly, $\{p_\varepsilon\}$ converges to a pressure field $p$ that also satisfies the properties stated. Thus,
the existence proof is completed, at least under the assumption $\bfv_*\in W^{1-\frac1q,q}(\partial\Omega)$. However, since this class of functions is dense in $W^{-\frac1q,q}(\partial\Omega)$ (e.g. \cite[Theorems 1, 2]{GSS}), the stated property follows from \eqref{2} and a standard density argument. \par\hfill$\square$
\setcounter{equation}{0}
\section{Existence and uniqueness of very weak solutions to the Navier-Stokes problem}
We are now interested in the same type of questions investigated in the previous section, but in the case of the following nonlinear problem
\be\ba{cc}\medskip\left.\ba{ll}\medskip
\Delta\bfv=\nabla p+\bfv\cdot\nabla\bfv+\Div\mathbb H\\ \medskip
\Div\bfv=0\ea\right\}\ \ \mbox{in $\Omega$\,,}\\
\bfv=\bfv_*\ \ \mbox{at $\partial\Omega$}\,, \ \ \Lim{|x|\to\infty}\bfv(x)=\0\,,
\ea
\eeq{3.1}
where $\mathbb H$ and $\bfv_*$ are given tensor and vector fields, respectively.
In analogy with \defref{1}, we give the following.
\Bd
The field $\bfv:\Omega\to \real^3$ is a very weak solution to \eqref{3.1} if for some $q\in (1,\infty)$ the following conditions are fulfilled.
\begin{itemize}
  \item [(a)] $\bfv\in L^q(\Omega_R)$, for all $R>\delta$\,;
  \item [(b)] $\bfv$ satisfies
 $$ 
\langle\bfv,\Delta\bfphi\rangle+\langle\bfv_*,\bfn\cdot{\nabla\bfphi}\rangle_{\partial\Omega}+\langle\bfv\cdot\nabla\bfphi,\bfv\rangle+\langle\mathbb H,\nabla\bfphi\rangle=0\,,\ \,\mbox{for all $\bfphi\in \calk_0^2(\Omega)$}\,;
$$
\item[(c)]\ $\langle\bfv,\nabla\zeta\rangle=0$\,,\ \ \mbox{for all $\zeta\in C_0^1(\real^3)$}\,;
\item[(d)]\ $\Liminf{|x|\to\infty}\int_{S^2}|\bfv(|x|,\omega)|{\rm d}\sigma_\omega=0\,.$  
\end{itemize}
\EDD{2}
The following result holds.
\Bt Let $q\in(6,\infty)$, and let $s$ and $\rho$ be as in \theoref{1}. Suppose that $\mathbb H$ and $\bfv_*$ satisfy
$$
\big(|x|^2\mathbb H\big)\in L^\infty(\Omega^s)\,;\ \ \mathbb H\in L^r(\Omega_{s})\,,\ \ \mbox{$\frac1r\le \frac1q+\frac13\,,\ r\in (3,\frac{q}2]$}\,;\ \ \bfv_*\in W^{-\frac1q,q}(\partial\Omega).
$$
There is $K=K(\Omega,\rho,s,q,r)>0$ such that if 
$$
\|\big(|x|^2\,\mathbb H\big)\|_{\infty,\Omega^s}+\|\mathbb H\|_{r,\Omega_{s}}+\|\bfv_*\|_{-\frac1q,q,\partial\Omega}\le K\,,
$$
then there exists a corresponding unique very weak solution $\bfv$ to \eqref{3.1} verifying
$$
(|x|\,\bfv)\in L^\infty(\Omega^s)\,,\ \bfv\in D^{1,q}(\Omega^\rho)\cap W^{1,r}(\Omega_{s,\rho})\,.
$$ 
Furthermore, we can find $p\in L^q(\Omega^\rho)\cap W^{-1,q}(\Omega_{2\rho})$ such that $(\bfv,p)$ satisfies \eqref{3.1}$_1$ in the sense of distributions. Finally, the pair $(\bfv,p)$ obeys the following estimate
\be
\ba{rl}\medskip
\|\big(|x|\,\bfv\big)\|_{\infty,\Omega^{s}}+\|\nabla\bfv\|_{q,\Omega^\rho}+\|\bfv\|_{1,r,\Omega_{s,\rho}}+&\!\!\!\!\|\bfv\|_{q,\Omega_{2\rho}}+\|p\|_{q,\Omega^\rho}+\|p\|_{-1,q,\Omega_\rho}\\
&\!\!\!\!
\le c\,\left(|\big(|x|^2\,\mathbb H\big)\|_{\infty,\Omega^s}+\|\mathbb H\|_{r,\Omega_{s}}+\|\bfv_*\|_{-\frac1q,q,\partial\Omega}\right)\,,
\ea 
\eeq{EST}
where $c=c(\Omega,\rho,s,q,r)>0$.
\ET{3}
{\em Proof.} Let $\hat{\bfv}$ be the very weak solution            constructed in \theoref{1} and corresponding to $\bfv_*$ and $\mathbb F\equiv\mathbb H$. Since $r<q$, from \eqref{2} we  have, in particular,
\be
\|(|x|\,\hat{\bfv})\|_{\infty,\Omega^{2\rho}}+\|\hat{\bfv}\|_{1,r,\Omega_{s,2\rho}}+\|\hat{\bfv}\|_{q,\Omega_{2\rho}}\le c\,{\sf D}\,,
\eeq{hat}
where
$$
{\sf D}:=\|\big(|x|^2\,\mathbb H\big)\|_{\infty,\Omega^s}+\|\mathbb H\|_{r,\Omega_{s}}+\|\bfv_*\|_{-\frac1q,q,\partial\Omega}\,.
$$
Taking into account that $r>3$, by embedding we deduce 
$$
\|\hat{\bfv}\|_{\infty,\Omega_{s,2\rho}}\le c\,\|\hat{\bfv}\|_{1,r,\Omega_{s,2\rho}}\big)\,.
$$
Thus, combining the latter with \eqref{hat}, we get 
\be
\|(|x|\,\hat{\bfv})\|_{\infty,\Omega^{s}}+\|\hat{\bfv}\|_{q,\Omega_{2\rho}}\le c\,{\sf D}\,.
\eeq{hat1}
Next, setting
\be
\bfu:=\bfv-\hat{\bfv}\,,
\eeq{3.2}
from \eqref{3.1} we formally obtain that $\bfu$ obeys  the nonlinear problem 
\be\ba{cc}\medskip\left.\ba{ll}\medskip
\Delta\bfu=\nabla p+\bfu\cdot\nabla\bfu+\hat{\bfv}\cdot\nabla\bfu+\bfu\cdot\nabla\hat{\bfv}\\ \medskip
\Div\bfu=0\ea\right\}\ \ \mbox{in $\Omega$\,,}\\
\bfu=\0\ \ \mbox{at $\partial\Omega$}\,.
\ea
\eeq{3.3}
We shall now show the existence of a very weak solution   to \eqref{3.3}, namely, a field $\bfu$ with
$$
(|x|\,\bfu)\in L^\infty(\Omega^\rho)\,,\ \bfu\in D^{1,q}(\Omega^\rho)\cap W^{1,r}(\Omega_{s,\rho})\cap L^q(\Omega_\rho)\,,
$$ 
satisfying 
\be
\ba{ll}\ms 
\langle\bfu,\Delta\bfphi\rangle+\langle\bfu\otimes\bfu+\hat{\bfv}\otimes\bfu+\bfu\otimes\hat{\bfv},\nabla\bfphi\rangle=0\,,\ \,\mbox{for all $\bfphi\in \calk_0^2(\Omega)$}\,,\\
\langle\bfu,\nabla\zeta\rangle=0\,,\ \ \mbox{for all $\zeta\in C_0^1(\real^3)$}\,.
\ea
\eeq{3.4}
This, in turn, by \eqref{2} and \theoref{1} implies that $\bfv$ is a very weak solution to \eqref{3.1} matching all the properties stated in the theorem.  
Consider the Banach space
$$
\cals^q=\cals^q(\Omega):= \{\bfu:\  \bfu\in L^q(\Omega_{\rho})\,,\ (|x|^2\bfu)\in L^\infty(\Omega^s);\ \bfu \
\mbox{satisfies \eqref{3.4}$_2$}\} 
$$ 
endowed with the norm
$$
\|\bfu\|_{\cals^q}:=\|\bfu\|_{q,\Omega_{\rho}}+\|(|x|^2\bfu)\|_{\infty,\Omega^s}\,,
$$
and, for a given $\epsilon_0>0$, its closed ball
$$
\cals_{\epsilon_0}^q=\cals_{\epsilon_0}^q(\Omega):=\{\bfu\in\cals^q:\ \|\bfu\|_{\cals^q}\le \epsilon_0\}\,.
$$
We next define the map
$$
{\sf M}:\bfw\in \cals^q_{\epsilon_0}\mapsto \bfu\in \cals^q
$$
where $\bfu$ satisfies 
\be\ba{ll}\medskip 
\langle\bfu,\Delta\bfphi\rangle+\langle\mathbb W,\nabla\bfphi\rangle=0\,,\ \,\mbox{for all $\bfphi\in \calk_0^2(\Omega)$}\,,
\\
\langle\bfu,\nabla\zeta\rangle=0\,,\ \ \mbox{for all $\zeta\in C_0^1(\real^3)$}
\ea\eeq{3.5}
with
\be
\mathbb W:=\bfw\otimes\bfw+\hat{\bfv}\otimes\bfw+\bfw\otimes\hat{\bfv}\,.
\eeq{3.6}
We will prove that, under the stated assumptions on the data, ${\sf M}$ has a fixed point. Let us begin to show that ${\sf M}$ is well defined. To this end, we observe that
$$\ba{rl}\medskip
\|(|x|^2\mathbb W)\|_{\infty,\Omega^s}+\|\mathbb W\|_{r,\Omega_s}\le &\!\!\!\!c\,\left[\|(|x|\bfw)\|_{\infty,\Omega^s}\big(\|(|x|\bfw)\|_{\infty,\Omega^s}+\|(|x|\hat{\bfv})\|_{\infty,\Omega^s}\big)\right.\\ 
&\hspace*{.7cm}\left.+\|\bfw\|_{2r,\Omega_s}\big(\|\bfw\|_{2r,\Omega_s}+\|\hat{\bfv}\|_{2r,\Omega_s}\big)\right]
\,.
\ea
$$
Since $2r\le q$, from  the latter and \eqref{hat1} we deduce
\be
\|(|x|^2\mathbb W)\|_{\infty,\Omega^s}+\|\mathbb W\|_{r,\Omega_s}\le c\|\bfw\|_{\cals^q}\left(\|\bfw\|_{\cals^q}+{\sf D}\right)\,.
\eeq{3.7}
Thus, $\mathbb W$ matches the assumptions of \theoref{1} and, consequently, \eqref{3.5} has a unique very weak solution $\bfu$ with an associated pressure field $p$ such that
$$
(|x|\,\bfu)\in L^\infty(\Omega^{2\rho})\,,\ \bfu\in D^{1,q}(\Omega^\rho)\cap W^{1,r}(\Omega_{s,\rho})\,,\ \ p\in L^q(\Omega^\rho)\cap W^{-1,q}(\Omega_\rho)\,,
$$ 
and satisfying
\be\ba{rl}\medskip
\|\big(|x|\,\bfu\big)\|_{\infty,\Omega^{2\rho}}+\|\nabla\bfu\|_{q,\Omega^\rho}+\|\bfu\|_{1,r,\Omega_{s,\rho}}+\|\bfu\|_{q,\Omega_{2\rho}}+&\!\!\!\!\|p\|_{q,\Omega^\rho}+\|p\|_{-1,q,\Omega_\rho}\\
&\!\!\!\!\le c\left(\|(|x|^2\mathbb W)\|_{\infty,\Omega^s}+\|\mathbb W\|_{r,\Omega_s}\right)\,,
\ea
\eeq{3.8}
Now, employing the assumption $r<q$ and \eqref{3.8}, we get
$$
\|\bfu\|_{1,r,\Omega_{s,2\rho}}\le\| \bfu\|_{1,r,\Omega_{s,\rho}}+\|\bfu\|_{q,\Omega_{2\rho}}+\|\nabla\bfu\|_{q,\Omega_{\rho,2\rho}} \le c\left(\|(|x|^2\mathbb W)\|_{\infty,\Omega^s}+\|\mathbb W\|_{r,\Omega_s}\right)\,,
$$
so that, since $r>3$, by embedding we deduce
$$
\|\bfu\|_{\infty,\Omega_{s,2\rho}} \le c\left(\|(|x|^2\mathbb W)\|_{\infty,\Omega^s}+\|\mathbb W\|_{r,\Omega_s}\right)\,,
$$
which,in turn, once combined with \eqref{3.8}, allows us to conclude
\be\ba{rl}\medskip
\|\big(|x|\,\bfu\big)\|_{\infty,\Omega^{s}}+\|\nabla\bfu\|_{q,\Omega^\rho}+\| \bfu\|_{1,r,\Omega_{s,\rho}}+\|\bfu\|_{q,\Omega_{2\rho}}\le c\left(\|(|x|^2\mathbb W)\|_{\infty,\Omega^s}+\|\mathbb W\|_{r,\Omega_s}\right)\,.
\ea
\eeq{3.9}
This proves, in particular, that the map ${\sf M}$ is well defined. We shall next show that ${\sf M}$ is a contraction.
From \eqref{3.7} and \eqref{3.9} it follows that, with a constant  $C=C(\Omega,\rho,s,q,r)$,
\be
\|\bfu\|_{\cals^q}\le C\,\|\bfw\|_{\cals^q}\left(\|\bfw\|_{\cals^q}+{\sf D}\right)\,,
\eeq{3.10}
which, in turn, by choosing $\bfw\in \cals_{\epsilon_0}^q$, entails
\be
\|\bfu\|_{\cals^q}\le C\,\big(\epsilon_0^2+{\sf D}\,\epsilon_0\big)\,.
\eeq{3.11}
As a result, if we take (for instance) $\epsilon_0$ such that
\be
{\sf D}=\frac{\epsilon_0}{2C}\,,\ \epsilon_0\le \frac1{4C+1}\,,
\eeq{3.12}
from \eqref{3.11} we infer that ${\sf M}$ is a self-map. It is then easy to show that ${\sf M}$ is also contractive. In fact, setting
$$
\bfu_i={\sf M}(\bfw_i)\,, \ i=1,2\,;\ \bfu:=\bfu_1-\bfu_2\,,\ \bfw:=\bfw_1-\bfw_2\,,
$$
from \eqref{3.5}--\eqref{3.6} we obtain
$$\ba{ll}\medskip 
\langle\bfu,\Delta\bfphi\rangle+\langle\hat{\mathbb W},\nabla\bfphi\rangle=0\,,\ \,\mbox{for all $\bfphi\in \calk_0^2(\Omega)$}\,,
\\
\langle\bfu,\nabla\zeta\rangle=0\,,\ \ \mbox{for all $\zeta\in C_0^1(\real^3)$}\,,
\ea$$
with
$$
\hat{\mathbb W}:=\bfw_1\otimes\bfw+\bfw\otimes\bfw_2+\hat{\bfv}\otimes\bfw+\bfw\otimes\hat{\bfv}\,.
$$
Therefore, proceeding as in \eqref{3.10}, we can show
$$
\|\bfu\|_{\cals^q}\le C\,\|\bfw\|_{\cals^q}\left(\|\bfw_1\|_{\cals^q}+\|\bfw_2\|_{\cals^q}+{\sf D}\right)\,,
$$
which, successively, by \eqref{3.12} implies
$$
\|\bfu\|_{\cals^q}\le C\,\|\bfw\|_{\cals^q}\big(2\epsilon_0+{\sf D}\big)\le C\,\epsilon_0(2+\frac1{2C})\|\bfw\|_{\cals^q}\le\frac12\,\|\bfw\|_{\cals^q}\,.
$$
The latter furnishes that ${\sf M}$ is contractive. This proves the existence of a unique solution to \eqref{3.4} in the ball $\cals^q_{\epsilon_0}(\Omega)$. However, from \eqref{3.7},\eqref{3.8}, \eqref{3.9}, \eqref{3.11} and the fact that $\bfu\equiv\bfw$, we show, on the one hand, that 
$$
\bfu\in D^{1,q}(\Omega^\rho)\cap W^{1,r}(\Omega_{s,\rho})\,,
$$
and, on the other hand, that $\bfu$ and $p$ satisfy the inequality
\be\ba{rl}\medskip
\|\big(|x|\,\bfu\big)\|_{\infty,\Omega^{s}}+\|\nabla\bfu\|_{q,\Omega^\rho}+\|\nabla \bfu\|_{r,\Omega_{s,\rho}}+\|\bfu\|_{q,\Omega_{2\rho}}+\|p\|_{q,\Omega^\rho}+\|p\|_{-1,q,\Omega_\rho}\le c\,{\sf D}\,.
\ea
\eeq{3.13}
Consequently, collecting all the properties proved for $\bfu,p$, and  taking into account \eqref{3.2} and \theoref{1}, we complete the proof of the theorem.\par\hfill$\square$\par
The next result provides an important consequence of \theoref{3}.
We recall the following notation for the Cauchy stress tensor:
$$
\mathbb T(\bfz,\phi):= 2\mathbb D(\bfz)-\phi\mathbb I\,,\ \ \mathbb D=\half(\nabla\bfz+(\nabla\bfz)^\top)\,,
$$
with $\mathbb I$ unit matrix.
\Bt Let $\Omega=\{x\in\real^3:\ |x|> R\}$, for some $R>0$. Let $q>6$ and suppose that $\bfv_*\in L^{\sf s}(\partial\Omega)$ for some ${\sf s}\ge \frac23\,q$. 
There is a positive constant $C_1=C_1(q,{\sf s})$ such that if 
\be
R\left(\int_{S^2}|\bfv_*(\omega)|^{\sf s} {\rm d}\sigma_\omega\right)^{\frac 1{\sf s}}\le  C_1\,,
\eeq{3.14}
then, problem \eqref{3.1}, with $\mathbb H\equiv \0$, has one and only one corresponding very weak solution, $\bfv$, of the type established in \theoref{3}. Moreover, $\bfv$ and the corresponding pressure field $p$ meet the following properties.
\begin{itemize}
\item[{\rm (i)}] $\bfv,p\in C^\infty(\Omega)$\,;
\item[{\rm (ii)}] There is a positive constant $C_2=C_2(q,{\sf s})$ such that
$$
\|(|x|\bfv)\|_{\infty,\Omega^{2R}}+\|\bfv\|_{3,\Omega_{6R}}\le C_2{R}\left(\int_{S^2}|\bfv_*(\omega)|^{\sf s} {\rm d}\sigma_\omega\right)^\frac1{\sf s}\,;
$$  
\item[{\rm (iii)}] For any $|\alpha|\ge0$, 
$$
D^\alpha\bfv(x)=O(|x|^{-1-|\alpha|})\,,\ \ D^\alpha p(x)=O(|x|^{-2-|\alpha|})\,,\ \ \mbox{as $|x|\to\infty$}\,; 
$$  
\item[{\rm (iv)}] If $\bfv_*\in W^{\frac32,2}(\partial\Omega)$,   then $(\bfv,p)\in W^{2,2}(\Omega_{\sf P})\times W^{1,2}(\Omega_{\sf P})$ for all ${\sf P>R}$ and obeys the  energy equality
$$ 
2\|\mathbb D(\bfv)\|_2^2 -\langle\bfv_*,\mathbb T(\bfv,p)\cdot\bfn\rangle_{\partial\Omega}+\half\langle |\bfv_*|^2,\bfv_*\cdot\bfn\rangle_{\partial\Omega}=0\,.
$$  
\end{itemize}
\ET{3.2}
{\em Proof.} Let $y:=x/R$, and define
$$
\bfv^\dagger(y):=R\,\bfv(Ry)\,,\ \ p^\dagger(y):=R^2\,p(Ry)\,, \ \ \bfv_*^\dagger(\omega):=R\,\bfv_*(\omega)\,,\ \ {\Omega}^\dagger:=\{y\in\real^3:\ |y|>1\}\,,
$$
where $\bfv,p$ and $\bfv_*$ satisfy \eqref{3.1} with $\mathbb H\equiv\0$.
Then, $\bfv^\dagger,p^\dagger$ and $\bfv_*^\dagger$ satisfy the following problem
\be
\ba{cc}\medskip\left.\ba{ll}\medskip
\Delta\bfv^\dagger=\nabla p^\dagger+\bfv^\dagger\cdot\nabla\bfv^\dagger\\ \medskip
\Div\bfv^\dagger=0\ea\right\}\ \ \mbox{in $\Omega^\dagger$\,,}\\
\bfv^\dagger(y)=\bfv_*^\dagger(y)\ \ \mbox{at $S^2$}\,,\ \ \Lim{|x|\to\infty}\bfv^\dagger(x)=\0\,,
\ea
\eeq{3.15}
where, of course, all the operators involved act on the $y$-variable.
We now apply \theoref{3} to \eqref{3.15}, with the choice 
$r=\frac q2$, $s=2$ and $\rho=3$. We then deduce that there is a constant $c_1=c_1(q)$ such that if, for some $q>6$,
\be
\|\bfv^\dagger_*\|_{-\frac1q,q,S^2}\le c_1\,,
\eeq{3.16}
then \eqref{3.15} has one and only one very weak solution with the properties stated in that theorem. In particular, \eqref{EST} entails
\be
\|(|y|\,\bfv^\dagger)\|_{\infty,\Omega^{\dagger2}}+\|\bfv^\dagger\|_{3,\Omega^\dagger_{6}}\le c_2\,\|\bfv_*^\dagger\|_{-\frac1q,q,S^2}\,,
\eeq{3.17}
where also the positive constant $c_2$ depends only on $q$.
Pick an arbitrary function  $\bfzeta\in W^{\frac1q,q'}(S^2)$, $\frac1{q'}+\frac1q=1$. Then, by classical results, we can find $\bfphi\in W^{1,q'}(\Omega^\dagger_{2})$, such that
\be
\bfphi|_{S^2}=\bfzeta\,,\ \ \bfphi|_{\partial B_2}=\0\,,\ \  \|\bfphi\|_{1,q',\Omega^\dagger_2}\le c_3\,\|\bfzeta\|_{\frac1q,q',S^2}\,,
\eeq{TR}
where $c_3=c_3(q)>0$.
We now have
\be
\left|\int_{S^2}\bfv_*^\dagger\cdot\bfzeta\right|\le \|\bfv_*^\dagger\|_{{\sf r},S^2}\|\bfzeta\|_{{\sf r}',S^2}\,,\ \ \mbox{$\frac{1}{\sf r}+\frac{1}{\sf r'}=1$}\,.
\eeq{TR1}
Since $q'<\frac65$, by well-known trace theorems (e.g. \cite[Theorem II.4.1]{Gab}) it follows that
$$
\|\bfzeta\|_{{\sf r}',S^2}\le c_4\,\|\bfphi\|_{1,q',\Omega^\dagger_{2}}\,,\ \mbox{for all ${\sf r}\in [\frac{2}3q,q]$}\,,
$$
with $c_4=c_4(q,{\sf r})>0$, which, once combined with \eqref{TR1} and  \eqref{TR}$_2$, furnishes
$$
\left|\int_{S^2}\bfv_*^\dagger\cdot\bfzeta\right|\le c_5\,\|\bfv_*^\dagger\|_{{\sf r},S^2}\|\bfzeta\|_{\frac1q,q',S^2}\,,\ \ \mbox{for all $\bfzeta\in W^{\frac1q,q'}(S^2)$}\,,
$$
where $c_5:=c_3\,c_4$. 
Therefore, we infer 
\be
\|\bfv_*^\dagger\|_{-\frac1q,q,S^2}\le c_5\|\bfv_*^\dagger\|_{{\sf r},S^2}\,, \ \mbox{for all ${\sf r}\in [\frac{2}3q,q]$}\,,
\eeq{3.18}
so that, a fortiori,
\eqref{3.16} is  satisfied if 
\be
\|\bfv_*^\dagger\|_{{\sf r},S^2}\le c_6\,,
\eeq{3.19}
with $c_6:=c_1/c_5$.  
Since
\be 
\|\bfv_*^\dagger\|_{{\sf r},S^2}=R\,\left(\int_{S^2}|\bfv_*(\omega)|^{\sf r}{\rm d}\sigma_\omega\right)^\frac1{\sf r}\,,
\eeq{3.20}
we deduce that a sufficient condition for the existence and uniqueness of solutions to the original problem in the stated function class, is given by
$$ 
R\,\left(\int_{S^2}|\bfv_*(\omega)|^{\sf r}{\rm d}\sigma_\omega\right)^\frac1{\sf r}\le c_6\,,\ \mbox{for some ${\sf r}\in [\frac{2}3q,q]$}\,.
$$
Clearly, by H\"older inequality, this restriction is certainly satisfied if $\bfv_*$ satisfies one of the type \eqref{3.14}, which thus completes the existence proof. We now pass to the proof of the properties (i)--(iv). We begin to notice that, since $\bfv\in L^3_{\rm loc}(\Omega)$, the validity of (i) follows  from \cite[Theorem IX.5.1]{Gab}. 
Next, from \eqref{3.17}, \eqref{3.18} and \eqref{3.20}, we also check at once the validity of (ii).
Furthermore, since $\bfv(x)=O(|x|^{-1})$ as $|x|\to\infty$, the properties stated in (ii) are a consequence of the proof of \cite[Lemma X.9.2]{Gab}. Finally, we observe that, for any fixed ${\sf P}>R$, $\bfv$ is a very weak solution to the following Stokes problem 
\be\ba{cc}\medskip\left.\ba{ll}\medskip
\Delta\bfv=\nabla p+\Div\mathbb S\\ \medskip
\Div\bfv=0\ea\right\}\ \ \mbox{in $\Omega_{\sf P}$\,,}\\
\bfv=\bfv_*\ \ \mbox{at $\partial\Omega_{\sf P}$}\,,\ \ \bfv=\bfv\ \ \mbox{at $\partial B_{\sf P}$}\,,
\ea
\eeq{3.21}
with $\mathbb S:=\bfv\otimes\bfv$. Because $\bfv\in L^q(\Omega_{\sf P})$, $q>6$, it follows that $\mathbb S\in L^s(\Omega_{\sf P})$, $s>3$. As a consequence, by the assumption on $\bfv_*$, property (i),  classical results on the Stokes problem \cite[Theorem IV.6.1 (b)]{Gab}, and, finally, the uniqueness of very weak solutions \cite[Theorem 3]{GSS} we deduce
$
\bfv\in W^{1,s}(\Omega_{\sf P})\,.
$
As a result, by embedding, we infer $\Div\mathbb S\in L^2(\Omega_{\sf P})$, so that, again by classical results for \eqref{3.21} \cite[Theorem IV.6.1 (a)]{Gab}, we get $(\bfv,p)\in W^{2,2}(\Omega_{\sf P})\times W^{1,2}(\Omega_{\sf P})$, for all ${\sf P}>R$. In view of all the above,  the asymptotic properties (iii), and recalling that 
\be
\Div\mathbb T(\bfz,\phi)=\Delta\bfz-\nabla \phi\,,\ \ \mbox{if $\Div\bfz=0$}\,,
\eeq{3.27}
the energy equality stated in (iv) is obtained by dot-multiplying by $\bfv$ both sides of \eqref{3.1} (with $\mathbb H=\0$ and $\Omega$ as in the statement of the theorem), integrating by parts over $\Omega_{\sf P}$, and then letting ${\sf P}\to\infty$. This completes te proof of the theorem.  
\par\hfill$\square$

\setcounter{equation}{0}
\section{On the asymptotic behavior of Leray solutions.}
We will now show how the existence theory developed in the previous section can provide insights into the behavior  at large spatial distances of Leray solutions in an exterior domain.\par
To this end, consider the following problem
\be
\ba{cc}\medskip\left.\ba{ll}\medskip
\Delta\bfu=\nabla {\sf p}+\bfu\cdot\nabla\bfu+\bfF\\ \medskip
\Div\bfu=0\ea\right\}\ \ \mbox{in $\Omega$\,,}\\
\bfu=\bfu_*\ \ \mbox{at $\partial\Omega$}\,,\ \ \Lim{|x|\to\infty}\bfu(x)=\0\,,
\ea
\eeq{4.1}
where $\bfF$ and $\bfu_*$ are assigned.
\Bd Let $\bfF\in D_0^{-1,2}(\Omega)$, $\bfu_*\in W^{\frac12,2}(\partial\Omega)$. The pair $(\bfu,{\sf p})$ is a {\em Leray solution} to \eqref{4.1} if the following conditions are satisfied.
\begin{itemize}
  \item [{\rm (a)}] $\bfu\in D^{1,2}(\Omega)\cap L^6(\Omega)$\,; \, \ $
{\sf p}\in L^2(\Omega_{R_1})\cap \big[L^3(\Omega^{R_2})+L^2(\Omega^{R_2})\big]\,,\, \, \mbox{arbitrary $R_1,R_2>\delta$}\,;
$
  \item [{\rm (b)}] $\bfu=\bfu_*$ in $W^{\frac12,2}(\partial\Omega)$\,; 
  \item [{\rm (c)}] $(\bfu,{\sf p})$ satisfies \eqref{4.1}$_{1,2}$ in the sense of distributions\,.
  \end{itemize}
\EDD{4.1}
\par
We recall the following result, proved in \cite[Theorem X.4.1 and Lemma X.1.1]{Gab}.
\Bp  There exists a constant $K=K(\Omega)>0$ such that if
$$
\Phi:=\left|\int_{\partial\Omega}\bfu_*\cdot\bfn\right|\le K
$$
there is at least one corresponding Leray solution $(\bfu,{\sf p})$ to \eqref{4.1}. Moreover,  $\bfu$ obeys the  {\em generalized energy inequality}:
\be
2\|\mathbb D(\bfu)\|_2^2- 2\langle \mathbb D(\bfu),\mathbb D(\bfU)\rangle+\langle\bfu\cdot\nabla\bfU,\bfu-\bfU\rangle+\langle \bfF,\bfu-\bfU\rangle\le 0\,,
\eeq{4.2}
where the field $\bfU$ is such that 
$$\bfU=\bfu_*\ \mbox{at $\partial\Omega$}\,;\,\ \Div\bfU=0\ \mbox{in $\Omega$}\,,
$$
and, for arbitrary $R>\delta$, 
$$
\bfU\in W^{1,2}(\Omega_R)\cap D^{1,q}(\Omega^R)\cap L^r(\Omega)\,,\ \ \mbox{for all $q\in (1,\infty)$ and all $r\in (\frac32,\infty)$}. 
$$
\EP{4.1}
\par
In order to investigate the asymptotic behavior of Leray solutions, we need a preliminary result. This will be established provided that the data have some further regularity. 
\Bl Suppose that 
\be
\bfF\in L^{\frac65}(\Omega)\cap L^2(\Omega_{\sf R})\,,\ \mbox{all ${\sf R}>\delta$}\,;\ \  \bfu_*\in W^{\frac32,2}(\partial\Omega)\,.
\eeq{4.3}
Then, Leray solutions satisfy 
\be 
(\bfu,{\sf p})\in W^{2,2}(\Omega_{\sf R})\times W^{1,2}(\Omega_{\sf R})\,,\ \ \mbox{for all ${\sf R}>\delta$}\,.
\eeq{4.5}
Furthermore, if $\bfu$ obeys \eqref{4.2}, then
for arbitrary $R>\delta$  the following inequality holds
\be
2\|\mathbb D(\bfu)\|_{2,\Omega^R}^2 -\langle\bfu,\mathbb T(\bfu,{\sf p})\cdot\bfn\rangle_{\partial B_R}+\half\langle |\bfu|^2,\bfu\cdot\bfn\rangle_{\partial B_R}+\langle \bfF,\bfu\rangle_{\Omega^R}\le0\,,
\eeq{4.4}
where $\bfn$ is the unit outer normal at $\partial\Omega^R$. 
\EL{4.1} 
{\em Proof.} We begin to notice that $L^\frac65(\Omega)\subset D_0^{-1,2}(\Omega)$. Moreover, since $\bfu\in L^6(\Omega)$, we find that, by H\"older inequality and \eqref{4.3}, the last term on the left-hand side of \eqref{4.4} is meaningful. Next, again by \eqref{4.3} and \cite[Theorem X.1.1]{Gab}, we infer the validity of \eqref{4.5},
which implies, in particular, that $(\bfu,{\sf p})$ satisfies \eqref{4.1}$_{1,2}$ almost everywhere in $\Omega$. Thus, by dot-multiplying both sides of \eqref{4.1} by $\bfU$, integrating by parts over $\Omega$ and recalling that $\bfU=\bfu_*$ at $\partial\Omega$, we get ($\langle\cdot,\cdot\rangle\equiv \langle\cdot,\cdot\rangle_\Omega$)
\be
\langle \bfu_*,\mathbb T(\bfu,{\sf p})\cdot\bfn\rangle_{\partial\Omega}-2\langle \mathbb D(\bfu),\mathbb D(\bfU)\rangle -\langle \bfu\cdot\nabla\bfu,\bfU\rangle-\langle\bfF,\bfU\rangle=0\,.
\eeq{4.6}
Notice that, in view of \eqref{4.5} and the summability properties of the extension $\bfU$, every term in \eqref{4.6} is well defined. We now observe that, also by an integration by parts that takes into account that $\bfu=\bfU=\bfu_*$ at $\partial\Omega$,
$$
\langle \bfu\cdot\nabla\bfu,\bfU\rangle=\langle \bfu\cdot\nabla(\bfu-\bfU),\bfU\rangle+\langle\bfu\cdot\nabla\bfU,\bfU\rangle=-\langle\bfu\cdot\nabla\bfU,\bfu-\bfU\rangle+\half \langle |\bfu_*|^2,\bfu_*\cdot\bfn\rangle_{\partial\Omega}\,.
$$
Therefore, from the latter and \eqref{4.6} we infer
\be
2\langle \mathbb D(\bfu),\mathbb D(\bfU)\rangle-\langle \bfu\cdot\nabla\bfU,\bfu-\bfU\rangle+\langle\bfF,\bfU\rangle=\langle \bfu_*,\mathbb T(\bfu,{\sf p})\cdot\bfn\rangle_{\partial\Omega}-\half \langle |\bfu_*|^2,\bfu_*\cdot\bfn\rangle_{\partial\Omega}
\eeq{4.7}
Employing \eqref{4.7} into \eqref{4.2}, we then conclude
\be
2\|\mathbb D(\bfu)\|_2^2-\langle \bfu_*,\mathbb T(\bfu,{\sf p})\cdot\bfn\rangle_{\partial\Omega}+\half \langle |\bfu_*|^2,\bfu_*\cdot\bfn\rangle_{\partial\Omega}+\langle\bfF,\bfu\rangle\le 0\,,
\eeq{4.8} 
which is the ``strong" form of the generalized energy inequality. Next, we dot-multiply both sides of \eqref{4.1}$_1$ by $\bfu$ and integrate by parts over the bounded domain $\Omega_R$ to get
\be\ba{rl}\medskip
-2\|\mathbb D(\bfu)\|_{2,\Omega_R}^2+\langle \bfu_*,\mathbb T(\bfu,{\sf p})\cdot\bfn\rangle_{\partial\Omega}+&\!\!\!\!\langle \bfu_*,\mathbb T(\bfu,{\sf p})\cdot\bfn\rangle_{\partial B_R}\\
&\!\!\!\!
-\half \langle |\bfu_*|^2,\bfu_*\cdot\bfn\rangle_{\partial\Omega}
-\half\langle |\bfu_*|^2,\bfu_*\cdot\bfn\rangle_{\partial B_R}
-\langle\bfF,\bfu\rangle_{\Omega_R}=0
\ea
\eeq{4.9}
where $\bfn$ is the unit outer normal at $\partial\Omega_R$. The inequality \eqref{4.4} then follows  by summing side-by-side \eqref{4.8} and \eqref{4.9}.
\par\hfill$\square$\par 
\par
We are now in a position to prove the main result of this section.
\Bt Let $\bfF$, $\bfu_*$ satisfy the assumptions of \lemmref{4.1}, with $\bfF$ of bounded support, and let $(\bfu,{\sf p})$ be a corresponding Leray solution obeying the generalized energy inequality \eqref{4.2}. Moreover, let ${\sf s}\in [\frac23q,\infty)$,  with $q>6$. Then, there exists a constant $C_0=C_0(q,{\sf s})>0$ such that if,   
\be
R\left(\int_{S^2}|\bfu(R,\omega)|^{\sf s}{\rm d}\sigma_\omega\right)^\frac1{\sf s}<{ C_0}\,,\, \ \ \mbox{for {\em some} sufficiently large $R$}\,,
\eeq{4.10} 
necessarily $(\bfu,{\sf p})$ must have the following behavior \mbox{as $|x|\to\infty$}:
$$
D^\alpha\bfu(x)=O(|x|^{-1-|\alpha|})\,,\ \ D^\alpha {\sf p}(x)=O(|x|^{-2-|\alpha|})\,,\,\ \mbox{arbitrary $|\alpha|\ge0$}\,. 
$$
\ET{4.1}
{\em Proof.} Let $R$ be fixed and so large that $\Omega_R$ contains the support of $\bfF$, and consider the  problem
\be
\ba{cc}\medskip\left.\ba{ll}\medskip
\Delta\bfu=\nabla {\sf p}+\bfu\cdot\nabla\bfu\\ \medskip
\Div\bfu=0\ea\right\}\ \ \mbox{in $\Omega^R$\,,}\\
\bfu=\bfv_*\ \ \mbox{at $\partial B_R$}\,,\ \ \Lim{|x|\to\infty}\bfu(x)=\0\,,
\ea
\eeq{4.11}
where we set $\bfv_*=\bfv_*(\omega):=\bfu(R,\omega)$. Notice that by  classical regularity results \cite[Theorem IX.5.1]{Gab} and \eqref{4.3}, it follows that $\bfu,{\sf p}\in C^\infty(\Omega^R)$ and $\bfv_*\in W^{\frac32,2}(\partial B_R)$. Moreover \cite[Theorem X.5.1]{Gab}
\be
\lim_{|x|\to\infty}|\bfu(x)|=0\,,\ \ \lim_{|x|\to\infty}|{\sf p}(x)|=0\,.
\eeq{4.12}
Let
$$
C_\star:=\min\{C_1,C_2\}
$$
where $C_i$, $i=1,2$, are the constants defined in \theoref{3.2}. We then infer that if
$$
R\left(\int_{S^2}|\bfv_*(\omega)|^{\sf s}{\rm d}\sigma_\omega\right)^{\frac1{\sf s}}\le C_\star\,,
$$
there exists a unique pair $(\bfv,p)$
such that
\be
\ba{cc}\medskip\left.\ba{ll}\medskip
\Delta\bfv=\nabla {p}+\bfv\cdot\nabla\bfv\\ \medskip
\Div\bfv=0\ea\right\}\ \ \mbox{in $\Omega^R$\,,}\\
\bfv=\bfv_*\ \ \mbox{at $\partial B_R$}\,,\ \ \Lim{|x|\to\infty}\bfv(x)=\0\,,
\ea
\eeq{4.13}
and satisfying all the properties stated in that theorem. In particular,
\be
\|(|x|\bfv)\|_{\infty,\Omega^{2R}}+\|\bfv\|_{3,\Omega_{6R}}\le C_\star{R}\left(\int_{S^2}|\bfv_*(\omega)|^{\sf s} {\rm d}\sigma_\omega\right)^\frac1{\sf s}\,.
\eeq{4.14} 
Dot-multiplying both sides of \eqref{4.13}$_1$ by $\bfu$, integrating by parts over $\Omega^{R,\rho}$ and recalling \eqref{3.27}, we get
\be
-2\langle \mathbb D(\bfu),\mathbb D(\bfv)\rangle_{\Omega^{R,\rho}}+\langle\bfv_*,\mathbb T(\bfv,p)\cdot\bfn\rangle_{\partial B_R}+\langle\bfu,\mathbb T(\bfv,p)\cdot\bfn\rangle_{\partial B_\rho}=\langle \bfv\cdot\nabla\bfv,\bfu\rangle_{\Omega^{R,\rho}}\,.
\eeq{4.15}
In view of the asymptotic properties of $\bfv$ and $p$ given in part (iii) of \theoref{3.2}, and \eqref{4.12}$_1$, we show
$$  
\lim_{\rho\to\infty}\langle\bfu,\mathbb T(\bfv,p)\cdot\bfn\rangle_{\partial B_\rho}=0\,.
$$
Moreover, from part (iv) of \theoref{3.2} and classical embedding theorems we deduce $\bfv\in L^\infty(\Omega_{\sf P})$, for all ${\sf P}>R$, which, combined with the asymptotic properties of $\bfv$ leads to $|x|\,|\bfv(x)|\le c$, $x\in\Omega^{R}$, where $c$ depends on $\bfv_*$. Therefore,
employing the classical  inequality \cite[Theorem II.6.1]{Gab}
\be
 \int_{\Omega^R}\frac{|\bfu(x)|^2}{|x|^2}\le 4\|\nabla\bfu\|_{2,\Omega^R}^2\,,
\eeq{4.16}
we get
$$
|\langle \bfv\cdot\nabla\bfv,\bfu\rangle_{\Omega^{R,\rho}}|\le \||\bfv|\bfu\|_{2,\Omega^R}\|\nabla\bfv\|_{2,\Omega^{R}}\le c\,
\|\nabla\bfv\|_{2,\Omega^{R}}\|\nabla\bfu\|_{2,\Omega^{R}}
\,.
$$
In view of the above, we pass to the limit $\rho\to\infty$ in \eqref{4.15} to find
\be
-2\langle \mathbb D(\bfu),\mathbb D(\bfv)\rangle_{\Omega^{R}}+\langle\bfv_*,\mathbb T(\bfv,p)\cdot\bfn\rangle_{\partial B_R}=\langle \bfv\cdot\nabla\bfv,\bfu\rangle_{\Omega^{R}}\,.
\eeq{4.17}
We now dot-multiply both sides of \eqref{4.11}$_1$ by $\bfv$ and integrate by parts over $\Omega^{R,\rho}$. We obtain
\be
-2\langle \mathbb D(\bfu),\mathbb D(\bfv)\rangle_{\Omega^{R,\rho}}+\langle\bfv_*,\mathbb T(\bfu,{\sf p})\cdot\bfn\rangle_{\partial B_R}+\langle\bfv,\mathbb T(\bfu,{\sf p})\cdot\bfn\rangle_{\partial B_\rho}=\langle \bfu\cdot\nabla\bfu,\bfv\rangle_{\Omega^{R,\rho}}\,.
\eeq{4.18}
Reasoning as before, we show
$$
|\langle \bfu\cdot\nabla\bfu,\bfv\rangle_{\Omega^{R,\rho}}|\le c\,
\|\nabla\bfv\|_{2,\Omega^{R}}\|\nabla\bfu\|_{2,\Omega^{R}}
\,.
$$
Furthermore, by (a) of \defref{4.1}, the regularity property of ${\sf p}$ and \eqref{4.12}$_2$ we know that
$$
\int_{\Omega^R}\left(|\nabla\bfu|^2+|{\sf p}|^3\right)<\infty\,,
$$
which implies the existence of a sequence $\{\rho_m\}$ such that
\be
\lim_{m\to\infty}\rho_m\int_{\partial B_{\rho_m}}\left(|\bfu|^6+|\nabla\bfu|^2+|{\sf p}|^3\right)=0\,.
\eeq{4.19}
By Schwarz and H\"older inequalities, and the asymptotic properties of $\bfv$, we infer
$$
|\langle\bfv,\mathbb T(\bfu,{\sf p})\cdot\bfn\rangle_{\partial B_{\rho_m}}|\le c \left(\|\nabla\bfu\|_{2,\partial B_{\rho_m}}+\rho_m^\frac13\,\|{\sf p}\|_{3,\partial B_{\rho_m}}\right)\,,
$$
so that, by \eqref{4.19}, we conclude
$$
\lim_{m\to\infty}\langle\bfv,\mathbb T(\bfu,{\sf p})\cdot\bfn\rangle_{\partial B_{\rho_m}}=0\,.
$$
Therefore, passing to the limit $m\to\infty$ in \eqref{4.18} with $\rho\equiv\rho_m$, it follows that
\be
-2\langle \mathbb D(\bfu),\mathbb D(\bfv)\rangle_{\Omega^{R}}+\langle\bfv_*,\mathbb T(\bfu,{\sf p})\cdot\bfn\rangle_{\partial B_R}=\langle \bfu\cdot\nabla\bfu,\bfv\rangle_{\Omega^{R}}\,.
\eeq{4.20}
Now, by property (iv) of \theoref{3.2},
\be
2\|\mathbb D(\bfv)\|_2^2 -\langle\bfv_*,\mathbb T(\bfv,p)\cdot\bfn\rangle_{\partial B_R}+\half\langle |\bfv_*|^2,\bfv_*\cdot\bfn\rangle_{\partial\Omega}=0\,.
\eeq{4.21}
whereas, by
\lemmref{4.1}, 
\be
2\|\mathbb D(\bfu)\|_{2,\Omega^R}^2 -\langle\bfv_*,\mathbb T(\bfu,{\sf p})\cdot\bfn\rangle_{\partial B_R}+\half\langle |\bfv_*|^2,\bfv_*\cdot\bfn\rangle_{\partial B_R}\le0\,,
\eeq{4.22}
Summing side-by-side \eqref{4.17}, \eqref{4.20}, \eqref{4.21} and \eqref{4.22} we deduce
\be
2\|\mathbb D(\bfw)\|_{2,\Omega^R}^2 +\langle |\bfv_*|^2,\bfv_*\cdot\bfn\rangle_{\partial B_R}\le \langle\bfv\cdot\nabla\bfv,\bfu\rangle_{\Omega_R} +\langle\bfu\cdot\nabla\bfu,\bfv\rangle_{\Omega_R}\,,
\eeq{4.23}
where $\bfw:=\bfv-\bfu$. We next observe that 
\be
\langle\bfv\cdot\nabla\bfv,\bfu\rangle_{\Omega_R} +\langle\bfu\cdot\nabla\bfu,\bfv\rangle_{\Omega_R}=\langle\bfw\cdot\nabla\bfv,\bfw\rangle_{\Omega^R}+\langle\bfw\cdot\nabla\bfv,\bfv\rangle_{\Omega^R}+\langle\bfu\cdot\nabla\bfv,\bfu\rangle_{\Omega^R}+\langle\bfu\cdot\nabla\bfu,\bfv\rangle_{\Omega^R}\,.
\eeq{4.24}
Integrating by parts, and using  $\bfw=\0$ at $\partial B_R$ we show
$$
\langle\bfw\cdot\nabla\bfv,\bfw\rangle_{\Omega^{R,\rho_m}}=-\langle\bfw\cdot\nabla\bfw,\bfv\rangle_{\Omega^{R,\rho_m}}+\langle\bfw\cdot\bfn,\bfv\cdot\bfw\rangle_{\partial B_{\rho_m}}\,.
$$
By H\"older inequality and the asymptotic properties of $\bfv$, we get 
$$
|\langle\bfw\cdot\bfn,\bfv\cdot\bfw\rangle_{\partial B_{\rho_m}}|\le \|\bfw\|^2_{6,\partial B_{\rho_m}}\|\bfv\|_{\frac32,\partial B_{\rho_m}}\le c\,\rho_m^\frac13\,\|\bfw\|^2_{6,\partial B_{\rho_m}}\,,
$$
so that passing to the limit $\rho_m\to\infty$, in view of \eqref{4.19}, we infer
\be 
\langle\bfw\cdot\nabla\bfv,\bfw\rangle_{\Omega^{R}}=-\langle\bfw\cdot\nabla\bfw,\bfv\rangle_{\Omega^{R}}\,.
\eeq{4.25}
In a similar fashion, we prove
\be
\langle\bfw\cdot\nabla\bfv,\bfv\rangle_{\Omega^R}=0\,,\ \ \langle\bfu\cdot\nabla\bfv,\bfu\rangle_{\Omega^R}=\langle|\bfv_*|^2,\bfv_*\cdot\bfn\rangle_{\partial B_R}-\langle\bfu\cdot\nabla\bfu,\bfv\rangle_{\Omega^R}\,.
\eeq{4.26}
Thus, taking into account \eqref{4.24}--\eqref{4.26},  \eqref{4.23} entails
\be
2\|\mathbb D(\bfw)\|_{2,\Omega_R}^2\le -\langle\bfw\cdot\nabla\bfw,\bfv\rangle_{\Omega^{R}}\,.
\eeq{4.27}
Further, by H\"older and Schwarz inequalities,
$$\ba{rl}\medskip
|\langle\bfw\cdot\nabla\bfw,\bfv\rangle_{\Omega^{R}}|\le |\langle\bfw\cdot\nabla\bfw,\bfv\rangle_{\Omega_{2R}}|&\!\!\!\!+|\langle\bfw\cdot\nabla\bfw,\bfv\rangle_{\Omega^{2R}}|
\\ &\!\!\!\!
\le \|\bfw\|_{6,\Omega^R}\|\nabla\bfw\|_{2,\Omega^R}\|\bfv\|_{3,\Omega_{2R}}+\||\bfv|\bfw\|_{2,\Omega^R}\|\nabla\bfw\|_{2,\Omega^R}
\ea
$$
As a result, using Sobolev inequality $\|\bfw\|_{6,\Omega^R}\le \gamma_S\|\nabla\bfw\|_{2,\Omega^R}$, with $\gamma_S$ numerical constant, and inequality \eqref{4.16} combined with \eqref{4.14}, we find that \eqref{4.27} leads to
$$ 
2\|\mathbb D(\bfw)\|_{2,\Omega^R}^2\le C_\star(2+\gamma_S)R\|\bfv_*\|_{{\sf s},S^2}\,\|\nabla\bfw\|_{2,\Omega^R}^2\,.
$$
Observing that $\|\nabla\bfw\|_2^2=2\|\mathbb D(\bfw)\|_2^2$, the theorem is then a consequence of the last inequality, with $C_0:= [C_\star(2+\gamma_S)]^{-1}$.\par\hfill$\square$\par
\Br It is worth remarking that, to date, Leray solutions are known to satisfy a property weaker than that requested by \eqref{4.10}. In fact, from the condition $\bfu\in L^6(\Omega)$  we deduce that there exists an unbounded sequence $\{R_m\}$ such that 
$$
\lim_{m\to\infty}R_m^\beta\left(\int_{S^2}|\bfu(R_m,\omega)|^{\sf s}{\rm d}\sigma_\omega\right)^{\frac1{\sf s}}=0\,,\ \ \mbox{for all}\, {\sf s}\in [1,6]\,,
$$
where $\beta=\frac12$. However, to satisfy \eqref{4.10} we would need $\beta=1$, for at least one ${\sf s}>4$.
\ER{4.1}
\ms\par
{\bf Acknowledgment.} Work  partially supported by NSF DMS Grant-2307811. The author would like to thank Mr. J.A.~Wein for helpful conversations.

\renewcommand{\theequation}{A.\arabic{equation}}\setcounter{equation}{0}
\section*{Appendix}
{\bf Lemma}. {\sl Let} $q\in (\frac32,\infty)$. {\sl Assume} ${\mathbb G}$, ${\bff}$, {\sl and} $g$ {\sl are 
 given second-order 
tensor, vector, and scalar fields, respectively, in} ${\Bbb R}^{3}${ \sl 
satisfying}
$$
\big(|x|^{2}\mathbb G\big)\in L^\infty(\real^3)\,;\  \bff\in W^{-1,q}(B_{2\rho})\,,\  g\in L^q(B_{2\rho})\,;\ \
\supp(\bff),\ \ \supp(g)\subset B_{\rho}
\,,
$$
{\sl Then, the problem}
\be\medskip\left.\ba{ll}\medskip
\Delta\bfu=\nabla {\sf p}+\Div\mathbb G+\bff\\ 
\Div\bfu=g\ea\right\}\ \ \mbox{\sl in $\real^3$\,,}
\eeq{A1}
{\sl has one and only one solution such that}
$$\ba{c}\medskip
\bfu\in D^{1,q}(\real^n), \ \ {\sf p}\in L^q(\real^n)\,, \, \
(|x|\,\bfu)\in L^{\infty}(B^{2\rho}).
\ea$$
{\sl Moreover, the solution satisfies the estimate}
\be
\|(|x|\,\bfu)\|_{\infty,B^{2\rho}}+\|\nabla\bfu\|_{q,\real^3}+\|{\sf p}\|
_{q,\real^3}
\le c\left(\|(|x|^2\,\mathbb G)\|_{\infty,\real^3}+\|\bff\|_{-1,q,B_{2\rho}}+
\|g\|_{q,B_\rho}\right)
\eeq{A.2}
{\sl with $c=c(q,\rho)$}.
\smallskip\par\noindent
{\em Proof.} A solution to \eqref{A1} is given by
$$
\bfu:=\bfu_1+\bfu_2+\bfu_3\,,\ \ {\sf p}=\tau_1+\tau_2+g\,,
$$
where
$$
\bfu_1=\mathbb U\star(\Div\mathbb G)\,,\ \bfu_2=\mathbb U\star\bff\,, \ \bfu_3=\nabla\cale\star g\,,\ \ \tau_1=\nabla\cale \star \Div\mathbb G\,,\ \ \tau_2=\nabla\cale \star \bff
$$
and $\mathbb U$, $\cale$ are Stokes and Laplace fundamental solutions, respectively, while $\star$ denotes convolution product. Once we show that $(\bfu,{\sf p})$ satisfies the stated properties, by known results on the Stokes problem it also follows that $(\bfu,{\sf p})$ is unique. By the Calder\`on-Zygmund theorem, we deduce ($\|\cdot\|_q\equiv\|\cdot\|_{q,\real^3}$)
\be
\|\nabla\bfu_1\|_{q}+\|\nabla\bfu_3\|_{q}+\|\tau_1\|_{q}\le c\,\big(\|\mathbb G\|_q+\|g\|_{q,B_\rho}\big)
\eeq{A.3}
Let $\bfzeta\in L^{q'}(\real^3)$ with $\supp (\bfzeta)\subset B_{\sf R}$. We then have
$$
\langle\nabla\bfu_2,\bfzeta\rangle= \langle \nabla\mathbb U\star\bff,\bfzeta\rangle=\langle \bff,\nabla\mathbb U\star\bfzeta\rangle\,.  
$$
Again by the Calder\`on-Zygmund theorem, we show 
$$
\|\nabla(\nabla\mathbb U\star\bfzeta)\|_{q'}\le c\,\|\bfzeta\|_{q'}\,,
$$
with $c$ independent of ${\sf R}$.
Since $q'\in(1,3)$, it follows that $\nabla\mathbb U\star\bfzeta\in D_0^{1,q'}(\real^3)$ \cite[Theorem II.7.6]{Gab}, 
so that, by the arbitrarity of ${\sf R}$ and $\bfzeta$, we conclude
\be
\|\nabla\bfu_2\|_q\le c\,\sup_{\mbox{\footnotesize $\bfphi\in D_0^{1,q'}(\real^3);\, \|\nabla\bfphi\|_{q'}=1$}}|\langle f,\bfphi\rangle|=:c\,\|\bff\|_{-1,q,\real^3}\,.
\eeq{A.4}
In a completely analogous way we infer
\be
\|\tau_2\|_q\le c\,\|\bff\|_{-1,q,\real^3}\,.
\eeq{A.5}
Let $\chi=\chi(|x|)$ be a smooth function  with compact support in $B_{2\rho}$ and equal to 1  on $\bar{B_{\frac43\rho}}$. Recalling that $q'<3$,  with the help of the Sobolev inequality we easily show for all 
$\bfphi \in D_{0}^{1,q'}({\Bbb R}^{3})$
\be
\parallel \chi\,\bfphi \parallel _{1,q^\prime ,B_{2\rho}}
\le c\,\mid \bfphi \mid _{1,q^\prime ,{\Bbb R}^3},
\eeq{A.6}
which implies, in particular, that $(\chi\,\bfphi)\in W_0^{1,q'}(B_{2\rho})$. Thus,
$$
\mid \langle{\bff},\bfphi\rangle\mid  = \mid \langle{\bff},\chi\,\bfphi\rangle
\mid \le \parallel {\bff}\parallel _{-1,q,B_{2\rho}}\parallel \chi\,\bfphi 
\parallel _{1,q^\prime ,B_{2\rho}}.
$$
and the latter along with \eqref{A.6}, yields
\be
\mid {\bff}\mid_{-1,q,{\Bbb R}^3}\le c\,\parallel {\bff}\parallel 
_{-1,q,B_{2\rho}}.
\eeq{A.7}
We next show the pointwise estimate. By the assumption on $\mathbb G$ and well known properties of $\mathbb U$,
$$
|\bfu_1(x)|=|\nabla\mathbb U\star\mathbb G|\le c\, \int_{\real^3}\frac{|\mathbb G(y)|}{|x- y|^2}{\rm d}y\le c\,\|(|x|^2\,\mathbb G)\|_{\infty,\real^3}\int_{\real^3}\frac{{\rm d}y}{|x-y|^2|y|^2}
$$
Therefore (e.g. \cite[Lemma II.9.2]{Gab}),
\be
\|(|x|\,\bfu_1)\|_{\infty,\real^3}\le  c\,\|(|x|^2\,\mathbb G)\|_{\infty,\real^3}\,.
\eeq{A.8}
Moreover, since for $|x|\ge 2\rho$ and $|y|\le \rho$ we find   $|x-y|\ge |x|-|y|\ge \half |x|$, taking into account that $|\nabla\cale(\xi)|\le |\xi|^{-2}$, $\xi\neq 0$,
we readily infer
\be
\|(|x|\,\bfu_3)\|_{\infty,B^{2\rho}}\le  c\,\|g\|_{q,B_\rho}\,,\ \ |x|\ge 2\rho\,.
\eeq{A.9}
Next, consider
\be
\mid {\bfu}_{2}(x)\mid x\mid \mid =\left\mid 
\int_{B_{\rho}}
{\bff}(y)
\cdot {\mathbb B}(x,y){\rm d}y\right\mid,
\eeq{V.7.16}
where, for $i,j=1,2,3$,
$$
(\mathbb B)_{ij}(x,y)=\chi(y)(\mathbb U)_{ij}(x-y)\mid x\mid .
$$
Since, as noted earlier on, for $y\in B_{\rho}$ and $\mid x\mid \ge 2\rho$ it is $|x-y|\ge\frac12|x|$, we get
$$
\mid {\mathbb U}(x-y)\mid +\mid \nabla {\mathbb U}(x-y)\mid \le c\mid x\mid 
^{-1}\,.
$$
which, combined with the properties of $\chi$, furnishes $\mathbb B(x,\cdot)\in W_0^{1,q'}(B_{2\rho})$ on the one hand and, on the other hand, $\|\mathbb B(x,\cdot)\|_{1,q',B_{2\rho}}\le c$, with the positive constant $c$ independent of $x\in B^{2\rho}$. From this and \eqref{V.7.16} we thus deduce
\be
\big| {\bfu}_{2}(x)\mid x\mid\big|\le\|\bff\|_{-1,q,B_{2\rho}}\|\mathbb B(x,\cdot)\|_{1,q',B_{2\rho}}
\le c\,\|\bff\|_{-1,q,B_{2\rho}}\,,\ \ |x|\ge 2\rho.
\eeq{A.12}
The estimate \eqref{A.2} then follows from \eqref{A.3}--\eqref{A.5}, \eqref{A.7}--\eqref{A.9} and \eqref{A.12}, which completes the proof.\par\hfill$\square$

\ed